\let\latexarabic\arabic
\let\latexdocument\document
\let\latexenddocument\enddocument

\RequirePackage[thmmarks]{ntheorem}
\makeatletter
\renewtheoremstyle{plain} 
{\item[\hskip\labelsep \theorem@headerfont ##1\ \textup{##2}\theorem@separator]} 
{\item[\hskip\labelsep \theorem@headerfont ##1\ \textup{##2}\ (##3)\theorem@separator]}
\makeatother

\documentclass[lineno]{biometrika}

\let\document\latexdocument
\let\enddocument\latexenddocument
\AtEndDocument{\printhistory}
\let\arabic\latexarabic
\def\rm{}

\usepackage{amsmath}

\usepackage{times}
\usepackage{bm}
\usepackage{natbib}

\usepackage[plain,noend]{algorithm2e}

\makeatletter
\renewcommand{\algocf@captiontext}[2]{\quad #1\algocf@typo. \AlCapFnt{}#2} 
\def\@algocf@capt@plain{top}
\renewcommand{\algocf@makecaption}[2]{%
  \addtolength{\hsize}{\algomargin}%
  \sbox\@tempboxa{\algocf@captiontext{#1}{#2}}%
  \ifdim\wd\@tempboxa >\hsize
    \hskip .5\algomargin%
    \parbox[t]{\hsize}{\algocf@captiontext{#1}{#2}}
  \else%
    \global\@minipagefalse%
    \hbox to\hsize{\box\@tempboxa}
  \fi%
  \addtolength{\hsize}{-\algomargin}%
}
\makeatother

\def\T{{ \mathrm{\scriptscriptstyle T} }}

	\usepackage{xcolor}
	\usepackage{stmaryrd}
	\usepackage{amsfonts}
	\usepackage{comment}
	\usepackage{subfig}

	\newcommand{\wh}[1]{\widehat{#1}}
	\newcommand{\wt}[1]{\widetilde{#1}}
	
	\newcommand{\ttinf}{2,\infty}
	\newcommand{\inftyinfty}{\infty,\infty}
	\newcommand{\op}{\operatorname{op}}
	\newcommand{\Frob}{\operatorname{F}}
	\newcommand{\diag}{\operatorname{diag}}
	
	\newcommand{\R}{\mathbb{R}}
	
	\newcommand{\Ex}{\operatorname{E}}
	
	\usepackage{color}
	\usepackage{hyperref}
	\hypersetup{
		colorlinks,
		citecolor=blue,
		linkcolor=purple,
		urlcolor=green
	}
	
	\usepackage[plain,noend]{algorithm2e}
		
	\usepackage[capitalise]{cleveref}

	\newcommand{\vmxobj}{\operatorname{Varimax}}
	\newcommand{\orthGroup}{\mathcal{O}}
	\newcommand{\signorthGroup}{\mathcal{P}}
	
	\newcommand{\ohp}{ \operatorname{o}_{\operatorname{P}}}
	\newcommand{\Ohp}{ \operatorname{O}_{\operatorname{P}}}

\begin{document}

\jname{Biometrika}
\jyear{}
\jvol{}
\jnum{}
\accessdate{}



\markboth{J. Cape}{On varimax asymptotics and spectral methods}

\title{On varimax asymptotics in network models and spectral methods for dimensionality reduction}

\author{J. Cape}
\affil{Department of Statistics, University of Wisconsin--Madison,\\ 1300 University Avenue, Madison, Wisconsin 53706 U.S.A. \email{jrcape@wisc.edu}}

\maketitle

\begin{abstract}
	Varimax factor rotations, while popular among practitioners in psychology and statistics since being introduced by H.~Kaiser, have historically been viewed with skepticism and suspicion by some theoreticians and mathematical statisticians. Now, work by K.~Rohe and M.~Zeng provides new, fundamental insight:~varimax rotations provably perform statistical estimation in certain classes of latent variable models when paired with spectral-based matrix truncations for dimensionality reduction. We build on this newfound understanding of varimax rotations by developing further connections to network analysis and spectral methods rooted in entrywise matrix perturbation analysis. Concretely, this paper establishes the asymptotic multivariate normality of vectors in varimax-transformed Euclidean point clouds that represent low-dimensional node embeddings in certain latent space random graph models. We address related concepts including network sparsity, data denoising, and the role of matrix rank in latent variable parameterizations. Collectively, these findings, at the confluence of classical and contemporary multivariate analysis, reinforce methodology and inference procedures grounded in matrix factorization-based techniques. Numerical examples illustrate our findings and supplement our discussion.
\end{abstract}



\section{Introduction}
\label{sec:introduction}


\subsection{Background}
\label{sec:background}

	Factor analysis has a long, rich history as a tenet of multivariate data analysis in statistics and psychology. Stemming from foundational work by F.~Galton, K.~Pearson, and C.~Spearman more than a century ago, factor analysis developed under the auspices of numerous researchers including L.~Thurstone, H.~Kaiser, H.~Rubin, M.G.~Kendall, and T.W.~Anderson \citep{vincent1953origin,bartholomew1995spearman}. In the twenty-first century, conducting factor analysis remains a mainstay of applied scientific domains, finding widespread use throughout psychology \citep{fabrigar1999evaluating,brown2015confirmatory}, economics \citep{chamberlain1983arbitrage,bai2008large}, and biology \citep{alter2000singular, argelaguet2018multi}.
	
	Much has been written about factor analysis. According to some:
	\begin{quote}
		Factor analysis is a branch of statistical theory concerned with the resolution of a set of descriptive variables in terms of a small number of categories or factors. \citep[p.~3]{holzinger1941factor}
	\end{quote}
	According to others:
	\begin{quote}
		Factor analysis is not a purely statistical technique; there is always a certain amount of guesswork in it which no fair-minded person would attempt to deny; but factor analysis is chiefly a research tool, and there is of necessity a lot of guesswork in all research [...] \citep[p.~117]{vincent1953origin}
	\end{quote}
	Ultimately, matrices lie at the heart of this apparent dichotomy. In their role as linear transformations, different matrices are used in practice to obtain differently rotated estimated factors or loadings, with the goal of obtaining scientifically meaningful, interpretable coordinate directions \citep{thurstone1935vectors}. Practitioners frequently select from a collection of diverse optimization criteria, leading to differently transformed point cloud configurations in low-dimensional Euclidean space via orthogonal or oblique matrices, with the most common criteria being promax, quartimax, equimax, and finally varimax \citep{kaiser1958varimax}, the focus of this paper.
	
	Historically, researchers have sought to develop theory that coherently justifies the aforementioned process of rotating estimates and the widespread subsequent empirical validation of exploratory analyses. Such efforts have largely centered around classical Gaussian data generating models and maximum likelihood-based procedures \citep{anderson1956statistical,anderson1962introduction}. In the case of varimax factor rotation matrices, which are orthogonal matrices, the process of rotating estimates neither alters the behavior of standard Gaussian likelihoods nor does it affect overall least-squares reconstruction error. Taken together, these observations have led some to believe that factor analysis is plagued by a fundamental conundrum.
	
	At long last, \cite{rohe2023vintage} sheds new light on the statistical properties of varimax factor rotations. The authors place their results in the context of certain so-called semiparametric factor models which produce data matrices in high-dimensional ambient space that possess low-dimensional population-level latent structure. Under certain identifiability conditions and assumptions that are markedly non-Gaussian, it is shown that pairing varimax rotations with the truncated singular value decomposition of an observed data matrix allows one to accurately estimate underlying generative variables. Put differently, \cite{rohe2023vintage} establishes that varimax rotations provably perform statistical estimation in certain classes of latent variable models when applied to spectral-based matrix truncations for dimensionality reduction.
	
	The novel theoretical insights provided in \cite{rohe2023vintage} have been enabled by the contemporary emergence and proliferation of entrywise matrix perturbation analysis \citep{abbe2020entrywise,cape2019two,cape2019signal,fan2018covariance,fan2021robust}, partially summarized in the monograph \cite{chen2021spectral}. At a high level, this line of work provides fine-grained high-resolution analysis of structured data matrices and their perturbations arising from statistical models. For such settings, matrix perturbations pertaining to estimation procedures can be precisely analyzed and significantly differ from deterministic worst-case behavior reflected by the classical theorems of C.~Davis, W.~Kahan, P.-\AA.~Wedin, and H.~Weyl found in matrix analysis textbooks.
	
	Random graph and network models, typified by the study of large-dimensional adjacency matrices that represent pairwise interactions between vertices, or nodes, via the presence or absence of edges, or links, are among the statistical settings where modern matrix perturbation tools are advancing methodology and inference capabilities. As addressed in the present paper, entrywise and indeed row-wise perturbation analysis can be used to obtain high-resolution uniform error rates and asymptotic fluctuations that hold at the level of individual node embeddings; such results improve upon coarser, average error rates available from classical matrix theory. Further still, certain random graph models that possess latent structure can be viewed as cousins of general factor models, thereby motivating concurrent interest in varimax factor rotations for obtaining interpretable coordinate representations in network settings.

\subsection{Contributions of this paper}
\label{sec:contributions}

	This paper studies the large-sample theoretical properties of varimax-transformed network embeddings in the context of latent variable stochastic graph models. We establish asymptotic multivariate normality results for low-dimensional node-specific embedding vectors which estimate generative latent variables. Our results hold for models of undirected networks, directed networks, and degree-corrected networks. Taken together, this paper builds upon the uniform error analysis of varimax conducted in \cite{rohe2023vintage}. Our results contribute to the contemporary proliferation of theory, methods, and applications for spectral graph clustering, and for multivariate analysis more generally. We provide further evidence that varimax factor rotations, when paired with spectral methods for dimensionality reduction, provide a coherent spectral estimation strategy. Spectral methods for dimensionality reduction offer more than merely exploratory data analysis capabilities:~they do, in fact, possess significance criteria and quantifiable estimation accuracy by virtue of our being able to discern the behavior of their asymptotic fluctuations. Moreover, given the simplicity, flexibility, and computational scalability afforded by truncated matrix factorizations, spectral methods serve as a useful alternative to likelihood-based methods in practice.
	
	This paper provides numerical examples but deliberately avoids real data analysis with varimax, since hundreds of existing papers have already demonstrated its empirical success and show post hoc improved interpretability of rotated estimated loadings and factors. Rather, this work is one of only a few existing papers to date that studies the theoretical properties of varimax factor rotations in statistical latent variable models. Interested readers can find further discussion of related twenty-first century real data applications and methodological advances in, for example, \cite{rohe2023vintage,chen2020new:arxiv:v1}.


\section{Preliminaries}
\label{sec:preliminaries}

	The generic factor model posits that observed data $X = [x_{1},\dots,x_{n}] \in \R^{p \times n}$ take the form
	\begin{equation}
		\label{eq:modern_factor_model}
		X = \mu 1_{n}^{\T} + L F^{\T} + E,
	\end{equation}
	where $\mu \in \R^{p}$ is an intercept vector, $1_{n} \in \R^{n}$ is a vector of all ones, $L = [l_{1},\dots,l_{p}]^{\T} \in \R^{p \times d}$ is a matrix of unobserved factor loadings, $F = [f_{1},\dots,f_{n}]^{\T} \in \R^{n \times d}$ is a matrix of unobserved individual factors, and $E \in \R^{p \times n}$ is a matrix of error terms having mean zero, typically presumed to be uncorrelated with or independent of $L$ and $F$. The problem of interest is to estimate $\mu$, $L$, and $F$, given knowledge of the true dimension $d \ll p, n$ or with an estimate thereof. Additional identifiability conditions are commonly assumed to enable inference in \cref{eq:modern_factor_model}, e.g.,~see \cite{anderson1956statistical,fan2021robust}.

	Classical matrix perturbation theory, based on Weyl's inequality for singular values \citep{weyl1912asymptotische} and Wedin's sine theta theorem for singular subspaces \citep{wedin1972perturbation}, establishes that the data matrix $X$, after centering, is well-approximated by $L F^{\T}$ when the error term $E$ has a comparatively small effect, for example, when its operator norm is small relative to that of $L F^{\T}$. Given that the matrix $L F^{\T}$ is necessarily low rank since $d \ll p, n$, a natural candidate estimator is to take a singular value truncation of $X$. Here, the challenge remains how to construct estimates $\wh{L}$ and $\wh{F}$ from $X$ so as to then quantify their statistical properties at a more refined level of granularity than mere operator norm or Frobenius norm bounds permit.
	
	In the above display, \cref{eq:modern_factor_model} is an example of a signal-plus-noise matrix model, in which an observed matrix $X$ additively decomposes into a signal term $\mu 1_{n}^{\T} + L F^{\T}$ and a noise term $E$. Spectral methods are particularly well suited for such models which naturally arise in numerous statistical applications including network analysis, image denoising, and covariance matrix estimation. In these settings, computing an appropriate truncated singular value decomposition often has the effect of denoising data, in an effort to approximate the underlying, unobserved signal. Our study of varimax factor rotations is rooted in this context.
	
	Under the semiparametric factor model of \cite{rohe2023vintage}, though here using slightly different notation, one observes a data matrix $A \in \R^{n \times m}$ consisting of independent entries such that, for latent factor matrices $Z \in \R^{n \times d}$, $Y \in \R^{m \times d}$ and a latent full rank matrix of coefficients $B \in \R^{d \times d}$, at the population level, in expectation,
	\begin{equation}
		\label{eq:semiparFM_expectation}
		\Ex(A \mid Z, Y) = Z B Y^{\T}.
	\end{equation}
	Upon observing $A$, the goal is to consistently estimate the matrices $Z$ and $Y$. Conditional on $Z$ and $Y$, writing $A = ZBY^{\T} + \left(A - ZBY^{\T}\right)$ illustrates signal-plus-noise matrix model form. Focusing here on $Z$ for the sake of discussion, \cite{rohe2023vintage} proposes and studies estimators of the form $\wh{Z} = n^{1/2} \wh{U} R_{\wh{U}}$, where the orthonormal matrix $\wh{U}$ consists of the top left singular vectors of $A$ and the orthogonal matrix $R_{\wh{U}}$ maximizes the varimax objective function
	\begin{equation}
		\label{eq:vmx_objective}
		\vmxobj\left(R; \wh{U}\right) = \sum_{j=1}^{d}\frac{1}{n}\sum_{i=1}^{n}\left[(\wh{U}R)_{ij}^{4} - \left\{\frac{1}{n}\sum_{k=1}^{n}(\wh{U}R)_{kj}^{2}\right\}^{2}\right]
	\end{equation}
	over the set of $d$-dimensional orthogonal matrices, denoted by $R \in \orthGroup(d) = \left\{R \in \R^{d \times d}: RR^{\T} = R^{\T}R = I_{d}\right\}$. Solutions $R_{\wh{U}}$, deemed varimax factor rotations, are by inspection only unique up to multiplication by signed permutation matrices, the latter of which are denoted by $\Pi \in \signorthGroup(d) = \left\{ \Pi \in \orthGroup(d): \Pi_{ij} \in \{-1,0,1\}\right\}$. Moreover, since $\wh{U}$ consists of orthonormal columns, maximizing \cref{eq:vmx_objective} is equivalent to maximizing the expression $\|\wh{U}R\|_{\ell_{4}}^{4}$ as a function of $R$, again over the set of $d \times d$ orthogonal matrices.

	For certain latent variable statistical matrix models, including but not limited to semi-sparse stochastic blockmodels to be specified in \cref{sec:main}, \cite{rohe2023vintage} establishes that the $n \times d$ matrix of latent factors $Z$ is well-estimated by $\wh{Z}$ in the sense that the uniform row-wise estimation error asymptotically vanishes with high probability as $n \rightarrow \infty$ when $m \asymp n$, i.e.,
	\begin{equation}
		\label{eq:Rohe_ttinf_SBM_summary}
		\underset{\Pi \in \mathcal{P}(d)}{\inf} \hspace{0.2em} \underset{1 \le i \le n}{\max} \|\wh{Z}_{i} - \Pi^{\T} \cdot Z_{i}\|_{\ell_{2}} = \ohp(1).
	\end{equation}
	In fact, more is true, since \cite{rohe2023vintage} provides explicit, context-specific big-O bounds that hold with high probability. To be clear, these bounds control the maximum Euclidean row norm, denoted alternatively at times by $\ell_{\ttinf}$ or $\|\cdot\|_{2,\infty}$ or $\|\cdot\|_{2\rightarrow\infty}$. Additional comments are provided in \cref{sec:coda}.
	
	Our point of departure is \cref{eq:Rohe_ttinf_SBM_summary}. We seek to quantify the asymptotic distribution and covariance structure of individual vectors in varimax-transformed Euclidean point clouds corresponding to spectral embeddings of networks. Our subsequent usage of the breve symbol to write $\breve{Z}_{i}$ in \cref{sec:main} plays the role of $Z_{i}$ in \cref{sec:preliminaries} provided the latter is assumed to be appropriately scaled for identifiability.

\section{Main results}
\label{sec:main}

\subsection{Setup}
\label{sec:setup}

	In an effort to streamline the readability of the main results, the following three subsections are devoted to stand-alone statements of three main theorems. Our asymptotic results hold for undirected, directed, and degree-corrected stochastic blockmodel graphs \citep{holland1983stochastic,karrer2011stochastic} as specified below. Throughout, latent cluster memberships are represented using the canonical set of $k$-dimensional standard basis vectors, denoted by $\{e_{1},\dots,e_{k}\}$.
	
	Concretely, this paper considers the following algorithm for obtaining low-dimensional varimax-rotated data embeddings from large observed data matrices.
	\begin{enumerate}
		\item[] \underline{Algorithm~1:~Dimensionality reduction for blockmodel-type data using varimax.}
		\item[--] Input:~An adjacency matrix $A \in \{0,1\}^{n \times m}$ and embedding dimension value $r \ll \min(n,m)$.
		\item Compute the best rank-$r$ approximation of $A$. Let $\wh{U}\in\R^{n \times r}$ denote the top-$r$ left singular vectors of $A$, and let $\wh{V}\in\R^{m \times r}$ denote the top-$r$ right singular vectors of $A$.
		\item Separately for $\wh{U}$ and $\wh{V}$, compute an orthogonal matrix that maximize the varimax objective function, i.e.,
		\begin{equation*}
			R_{\wh{U}} \in \underset{R \in \orthGroup(r)}{\arg\max}\; \vmxobj(R; \wh{U}), \qquad
			R_{\wh{V}} \in \underset{R \in \orthGroup(r)}{\arg\max}\; \vmxobj(R; \wh{V}).
		\end{equation*}
		\item[--] Output:~The estimated latent factor matrices $\wh{Z} = n^{1/2}\wh{U}R_{\wh{U}}$ and $\wh{Y} = m^{1/2}\wh{V}R_{\wh{V}}$.
		\item[] (Optional:~The value $r$ could instead be selected in a data-driven manner, for example, based on inspecting the singular values of $A$.)
	\end{enumerate}
	In this paper, $\wh{Z}$ and $\wh{Y}$ are viewed as statistical estimators, and their row-wise distributional properties are quantified in the large-data limit. In practice, after obtaining $\wh{Z}$ and $\wh{Y}$, analysts might choose to visualize them using pair-pair plots, to apply a clustering algorithm to their rows, or to use them as plug-in estimates, depending on the problem at hand.

	The above algorithm is sufficient for the purposes of this paper. We emphasize that Algorithm~1 represents a special case of the more general \texttt{vsp} algorithm developed in \cite{rohe2023vintage}. Therein, the authors provide supporting discussion and guidance regarding the utility of optional centering, recentering, and scaling steps that are appropriate for more general data settings but are beyond the scope of the present paper.
	
\subsection{Undirected stochastic blockmodels}
\label{sec:main_SBM}

	Fix $k \ge 2$. Let $\pi$ be a $k$-dimensional probability vector with positive entries summing to unity. Let $B \in (0,1) ^{k \times k}$ be a full-rank symmetric matrix. Given any positive integer $n$, for all $i \in \{1,\dots,n\}$ generate independent, identically distributed $k$-dimensional latent indicator vectors $Z_{i} \sim \mathcal{G} = \operatorname{Multinomial}(1; \pi)$, arranged in the matrix $Z = [Z_{1},\dots,Z_{n}]^{\T} \in \{0,1\}^{n \times k}$.
	
	Conditional on $Z$, sample a symmetric adjacency matrix $A \in \{0,1\}^{n \times n}$ whose upper triangular entries are independent Bernoulli random variables with edge probabilities given by $\Ex(A \mid Z) = \rho Z B Z^{\T}$, where $\rho = \rho_{n} \in (0,1]$ denotes a sparsity term. Define $\wh{Z}_{i}$ as the $i$-th row of $\wh{Z} = n^{1/2}\wh{U}R_{\wh{U}}$ obtained from Algorithm~1. Let $\breve{Z}_{i} = \diag(\pi)^{-1/2} Z_{i}$.

\begin{theorem}
	\label{thrm:main_SBM}
	Assume the undirected stochastic blockmodel setting of \cref{sec:main_SBM} with $n \rho = \omega(\log^{c} n)$ for some sufficiently large constant $c > 1$ where either $\rho = 1$ or $\rho \rightarrow 0$ with limit $\rho_{\infty}$. It suffices to set $c=20$. There exists a sequence of signed permutation matrices $(\Pi_{Z}) \subset \signorthGroup(k)$ such that for any choice of fixed index $i$, as $n \rightarrow \infty$, the sequence of random vectors of the form
	\begin{equation}
	\label{eq:main_SBM_Zhat_CLT}
		(n\rho)^{1/2}\left\{\Pi_{Z} \cdot \wh{Z}_{i} - \breve{Z}_{i} \mid Z_{i} = e_{\ell} \right\}
	\end{equation}
	converges in distribution to a multivariate Gaussian random vector with mean zero and covariance matrix given by
	\begin{equation}
		\label{eq:thrm_covariance_SBM_Zhat}
		\left\{\diag(\pi)^{-1/2} \cdot T_{B}^{-1} \cdot J \cdot \Delta_{X}^{-1}\right\} \cdot \Ex_{\xi}\{ g_{\rho_{\infty}}(\xi, T_{B} e_{\ell}) \} \cdot \left\{\Delta_{X}^{-1} \cdot J \cdot (T_{B}^{-1})^{\T} \cdot \diag(\pi)^{-1/2}\right\}.
	\end{equation}
	In \cref{eq:thrm_covariance_SBM_Zhat}, $T_{B} = |\Lambda_{B}|^{1/2} U_{B}^{\T}$ is defined in terms of the spectral decomposition of $B$ such that $B = T_{B}^{\T} J T_{B}$ where $J = \diag(1_{p},-1_{q}) \in \R^{k \times k}$ when $B$ has $p$ positive eigenvalues and $q$ negative eigenvalues. Furthermore, $\Delta_{X} = \Ex(\xi\xi^{\T})$ where $\xi \sim T_{B} \mathcal{G}$, and $g_{\rho_{\infty}}(\xi, x) = (x^{\T}J\xi) \cdot (1-\rho_{\infty}x^{\T}J\xi) \cdot \xi\xi^{\T}$ where $\rho_{\infty} \in \{0,1\}$.
\end{theorem}

\subsection{Directed stochastic blockmodels}
\label{sec:main_DirSBM}

		Fix $k \ge 2$. Let $\pi_{Z}$ and $\pi_{Y}$ be two $k$-dimensional probability vectors each having positive entries summing to unity. Let $B \in (0,1)^{k \times k}$ be a full-rank matrix though not necessarily symmetric. Given any positive integer $n$, for all $i \in \{1,\dots,n\}$ generate independent, identically distributed $k$-dimensional latent indicator vectors $Z_{i} \sim \mathcal{G} = \operatorname{Multinomial}(1; \pi_{Z})$, arranged in the matrix $Z = [Z_{1},\dots,Z_{n}]^{\T} \in \{0,1\}^{n \times k}$. Similarly, for all $j \in \{1,\dots,n\}$, generate independent, identically distributed $k$-dimensional latent indicator vectors $Y_{j} \sim \mathcal{H} = \operatorname{Multinomial}(1; \pi_{Y})$, arranged in the matrix $Y = [Y_{1},\dots,Y_{n}]^{\T} \in \{0,1\}^{n \times k}$.
		
		Conditional on $Z$ and $Y$, sample a directed adjacency matrix $A \in \{0,1\}^{n \times n}$ consisting of independent Bernoulli entries with edge probabilities given by $\Ex(A \mid Z, Y) = \rho Z B Y^{\T}$, where $\rho = \rho_{n} \in (0,1]$ denotes a sparsity term. Define $\wh{Z}_{i}$ as the $i$-th row of $\wh{Z} = n^{1/2}\wh{U}R_{\wh{U}}$ obtained from Algorithm~1. Similarly, define $\wh{Y}_{j}$ as the $j$-th row of $\wh{Y} = n^{1/2}\wh{V}R_{\wh{V}}$ obtained from Algorithm~1. Let $\breve{Z}_{i} = \diag(\pi_{Z})^{-1/2} Z_{i}$ and $\breve{Y}_{j} = \diag(\pi_{Y})^{-1/2} Y_{j}$.

\begin{theorem}
	\label{thrm:main_DirSBM}
	Assume the directed stochastic blockmodel setting of \cref{sec:main_DirSBM} with $n \rho = \omega(\log^{c} n)$ for some sufficiently large constant $c > 1$ where either $\rho = 1$ or $\rho \rightarrow 0$ with limit $\rho_{\infty}$. It suffices to set $c=20$. There exists a sequence of signed permutation matrices $(\Pi_{Z}) \subset \signorthGroup(k)$ such that for any choice of fixed index $i$, as $n \rightarrow \infty$, the sequence of random vectors of the form
	\begin{equation}
	\label{eq:main_DirSBM_Zhat_CLT}
		(n\rho)^{1/2}\left\{\Pi_{Z} \cdot \wh{Z}_{i} - \breve{Z}_{i} \mid Z_{i} = e_{\ell} \right\}
	\end{equation}
	converges in distribution to a multivariate Gaussian random vector with mean zero and covariance matrix given by
	\begin{equation}
		\label{eq:thrm_covariance_DirSBM_Zhat}
		\left\{\diag(\pi_{Z})^{-1/2} \cdot T_{\Xi}^{-1} \cdot \Delta_{\Upsilon}^{-1}\right\} \cdot \operatorname{E}_{\upsilon}\left\{ g_{\rho_{\infty}}(\upsilon, T_{\Xi} e_{\ell}) \right\} \cdot \left\{\Delta_{\Upsilon}^{-1} \cdot (T_{\Xi}^{-1})^{\T} \cdot \diag(\pi_{Z})^{-1/2}\right\}.
	\end{equation}
	
	Simultaneously, there exists a sequence of signed permutation matrices $(\Pi_{Y}) \subset \signorthGroup(k)$ such that for any choice of fixed index $j$, as $n \rightarrow \infty$, the sequence of random vectors of the form
	\begin{equation}
	\label{eq:main_DirSBM_Yhat_CLT}
		(n\rho)^{1/2}\left\{\Pi_{Y} \cdot \wh{Y}_{j} - \breve{Y}_{j} \mid Y_{j} = e_{\ell} \right\}
	\end{equation}
	converges in distribution to a multivariate Gaussian random vector with mean zero and covariance matrix given by
	\begin{equation}
		\label{eq:thrm_covariance_DirSBM_Yhat}
		\left\{\diag(\pi_{Y})^{-1/2} \cdot T_{\Upsilon}^{-1} \cdot \Delta_{\Xi}^{-1}\right\} \cdot \Ex_{\xi}\left\{ g_{\rho_{\infty}}(\xi, T_{\Upsilon} e_{\ell}) \right\} \cdot \left\{\Delta_{\Xi}^{-1} \cdot (T_{\Upsilon}^{-1})^{\T} \cdot \diag(\pi_{Y})^{-1/2}\right\}.
	\end{equation}
	In \cref{eq:thrm_covariance_DirSBM_Zhat,eq:thrm_covariance_DirSBM_Yhat}, the matrices of coefficients $T_{\Xi} = S_{B}^{1/2}U_{B}^{\T}$ and $T_{\Upsilon} = S_{B}^{1/2}V_{B}^{\T}$ are defined in terms of the singular value decomposition of $B$, written $B = U_{B} S_{B} V_{B}^{\T}$. Here, $\Delta_{\Xi} = \Ex(\xi\xi^{\T})$ where $\xi \sim T_{\Xi}\mathcal{G}$, $\Delta_{\Upsilon} = \Ex(\upsilon\upsilon^{\T})$ where $\upsilon \sim T_{\Upsilon}\mathcal{H}$, and $g_{\rho_{\infty}}(\alpha, \beta) = (\beta^{\T}\alpha) \cdot (1-\rho_{\infty}\beta^{\T}\alpha) \cdot \alpha\alpha^{\T}$ with $\rho_{\infty} \in \{0,1\}$ takes $k$-dimensional input vectors $\alpha$ and $\beta$.
\end{theorem}
	
	In general, the asymptotic covariance matrices differ for \cref{eq:main_DirSBM_Zhat_CLT} and \cref{eq:main_DirSBM_Yhat_CLT}. The asymptotic distribution of \cref{eq:main_DirSBM_Zhat_CLT} depends on the distribution $\mathcal{H}$, and similarly, \cref{eq:main_DirSBM_Yhat_CLT} depends on the distribution $\mathcal{G}$, manifest in the covariance expressions \cref{eq:thrm_covariance_DirSBM_Zhat,eq:thrm_covariance_DirSBM_Yhat}.

	For simplicity, \cref{thrm:main_DirSBM} is stated and proved for large, square adjacency matrices. It can be shown that asymptotic normality continues to hold for rectangular adjacency matrices $A \in \{0,1\}^{n \times m}$ in the large-system asymptotic regime $n, m \rightarrow \infty$ when $n$ and $m$ grow at the same rate, requiring only minor modifications to the current notation and proof.

\subsection{Degree-corrected stochastic blockmodels}
\label{sec:main_DCSBM}

		Fix $k \ge 2$. Let $\pi$ be a $k$-dimensional probability vector with positive entries summing to unity. Let $B \in (0,1)^{k \times k}$ be a full-rank symmetric matrix. Given any positive integer $n$, for all $i \in \{1,\dots,n\}$ generate independent, identically distributed $k$-dimensional latent indicator vectors $Z_{i} \sim \mathcal{G} = \operatorname{Multinomial}(1; \pi)$, arranged in the matrix $Z = [Z_{1},\dots,Z_{n}]^{\T} \in \{0,1\}^{n \times k}$. Independently of $Z$, generate independent, identically distributed scalar degree parameters $\theta_{i} \sim \mathcal{H}$ from a distribution $\mathcal{H}$ with compact support contained in the open unit interval, to form the vector of degree-heterogeneity parameters $\theta = (\theta_{1},\dots,\theta_{n})^{\T}$.
		
		Conditional on $Z$ and $\theta$, sample a symmetric adjacency matrix $A \in \{0,1\}^{n \times n}$ whose upper triangular entries are independent Bernoulli random variables with edge probabilities given by $\Ex(A \mid Z, \theta) = \rho \diag(\theta) Z B Z^{\T} \diag(\theta)$, where $\rho = \rho_{n} \in (0,1]$ denotes a sparsity term. Define $\wh{Z}_{i}$ as the $i$-th row of $\wh{Z} = n^{1/2}\wh{U}R_{\wh{U}}$ obtained from Algorithm~1. Let $\breve{Z}_{i}^{\prime} = \diag(\eta)^{-1/2} \cdot (\theta_{i} Z_{i})$ where $\eta = \Ex(\theta_{1}^{2}Z_{1}^{2}) \in\R^{k}$.
		
	\begin{theorem}
		\label{thrm:main_DCSBM}
		Assume the degree-corrected stochastic blockmodel setting of \cref{sec:main_DCSBM} with $n \rho = \omega(\log^{c} n)$ for some sufficiently large constant $c > 1$ where either $\rho = 1$ or $\rho \rightarrow 0$ with limit $\rho_{\infty}$. It suffices to set $c=20$. There exists a sequence of signed permutation matrices $(\Pi_{Z}) \subset \signorthGroup(k)$ such that, for any choice of fixed index $i$ and any vector $a \in \R^{k}$, it holds that
		\begin{equation}
			\label{eq:main_DCSBM_Zhat_CLT}
			\operatorname{pr}\left\{(n\rho)^{1/2}\left(\Pi_{Z} \cdot \wh{Z}_{i} - \breve{Z}_{i}^{\prime}\right) \le a \right\}
			\longrightarrow \int_{\operatorname{supp}\mathcal{F}} \Phi\left\{a, \Gamma_{\rho_{\infty}}(b)\right\}\operatorname{d}\mathcal{F}(b)
		\end{equation}
		in the large-$n$ limit, where $\mathcal{F}$ denotes the product distribution $\mathcal{G}\mathcal{H}$ and $\Phi\left\{\cdot,\cdot\right\}$ denotes the mean zero multivariate Gaussian cumulative distribution function evaluated at the vector $a \in \R^{k}$ with covariance matrix $\Gamma_{\rho_{\infty}}(b)$, for $b \in \operatorname{supp}\mathcal{F}$, given by
		\begin{equation}
			\left\{\diag(\eta)^{-1/2} \cdot T_{B}^{-1} \cdot J \cdot \Delta_{X^{\prime}}^{-1}\right\} \cdot \Ex_{\xi^{\prime}}\{ g_{\rho_{\infty}}(\xi^{\prime}, b)\} \cdot \left\{\Delta_{X^{\prime}}^{-1} \cdot J \cdot (T_{B}^{-1})^{\T} \cdot \diag(\eta)^{-1/2}\right\}.
		\end{equation}
		Above, $T_{B} = |\Lambda_{B}|^{1/2} U_{B}^{\T}$ satisfies $B = T_{B}^{\T} J T_{B}$ where $J = \diag(1_{p},-1_{q}) \in \R^{k \times k}$ when $B$ has $p$ positive eigenvalues and $q$ negative eigenvalues. Further, $\Delta_{X^{\prime}} = \Ex(\xi^{\prime}\xi^{\prime\T})$ where $\xi^{\prime} \sim T_{B} \mathcal{G}\mathcal{H}$, and $g_{\rho_{\infty}}(\xi^{\prime}, b) = (b^{\T}J\xi^{\prime}) \cdot (1-\rho_{\infty}b^{\T}J\xi^{\prime}) \cdot \xi^{\prime}\xi^{\prime\T}$ where $\rho_{\infty} \in \{0,1\}$.
	\end{theorem}
	
\subsection{Coda to main results}
\label{sec:coda}

	Several remarks are in order to further contextualize and clarify the main results. First, regarding notation, the literature on stochastic blockmodels commonly writes $Z_{i}$ to denote the unscaled community membership vector for node $i$ taking values in the set of $k$-dimensional standard basis vectors. Above, our notation $\breve{Z}_{i}$ indicates that the community membership vector is scaled for the sake of identifiability, though for simplicity we drop the breve symbol in some of the surrounding discussion here. Second, writing the entrywise matrix statement $\wh{Z} \approx Z\Pi$ corresponds to writing the column vector statement $\wh{Z}_{i} \approx \Pi^{\T}Z_{i}$ and hence $\Pi \wh{Z}_{i} \approx Z_{i}$ uniformly for all vector indices $i \in \{1,\dots,n\}$.
	
	In the main theorems, each sequence of signed permutation matrices $(\Pi) \subset \signorthGroup(k)$ is implicitly indexed by $n$ and does not depend on the choice of node index $i$. The need to account for these matrices is unavoidable in general and discussed further in \cref{sec:discussion}.

	The big-O rate established in \cite{rohe2023vintage} for stochastic blockmodels, which implies the uniform error bound in \cref{eq:Rohe_ttinf_SBM_summary}, is approximately $\Ohp\left\{(n\rho)^{-1/4}\right\}$, ignoring logarithmic terms. In contrast, the results in this paper suggest the non-uniform row-wise $\ell_{2}$ norm bound $\Ohp\left\{(n\rho)^{-1/2}\right\}$, ignoring logarithmic terms, noting that asymptotic normality is driven by the stabilization of a linear combination of independent latent bounded random variables. To be clear, this paper obtains asymptotic normality for individual low-dimensional embedding vectors, thereby building on the results in \cite{rohe2023vintage}, but without obtaining an improved uniform error bound for \cref{eq:Rohe_ttinf_SBM_summary}.
	
	The independence assumption between $Z_{i}$ and $Y_{i}$ in each tuple $(Z_{i}, Y_{i})$ in \cref{sec:main_DirSBM} is conventional, though it is perhaps possible that versions of \cref{eq:main_DirSBM_Zhat_CLT,eq:main_DirSBM_Yhat_CLT} continue to hold under certain forms of dependence. The choice of constant $c=20$ in the growth condition $n\rho = \omega(\log^{c} n)$, i.e.,~$n\rho/(\log^{c}n) \rightarrow \infty$ where $\log^{c} n = (\log n)^{c}$, is made for convenience and is based on the perturbation analysis collectively found in \cite{cape2019signal,tang2022asymptotically,rohe2023vintage}. The asymptotic covariance matrices depend on the graph sparsity regime $\rho_{\infty} \in \{0,1\}$ and convey larger asymptotic variability when $\rho_{\infty} = 0$ compared to the dense graph regime $\rho_{\infty} = 1$.
	
	A key challenge in going from the $\ell_{\ttinf}$ bound in \cite{rohe2023vintage} to a row-wise distributional guarantee is to decompose the matrix difference between the sample quantity and the population-level quantity, namely the difference between the estimated latent factors and ground-truth latent factors, into tractable dominant and residual terms that reflect the underlying signal-plus-noise structure. Beyond simply identifying the dominant leading-order behavior in the perturbation analysis, it must then be appropriately oriented and decomposed so as to obtain explicit formulas for the asymptotic mean vector and covariance matrix in terms of the data-generating model parameters. Non-trivial bookkeeping and analysis involving various orthogonal transformations is required. Further details and discussions can be found in the proofs which are provided in the Supplementary Material.
	
\section{Numerical examples}
\label{sec:numerics}

\subsection{Four-block undirected stochastic blockmodel graphs}
\label{sec:numerics_undirected}

	Here, we simulate $100$ undirected stochastic blockmodel graphs as defined in \cref{sec:main_SBM}, each having $n = 5000$ nodes, sparsity term $\rho=1$, and with block connectivity matrix and membership probability vector given by
	\begingroup
	\renewcommand*{\arraystretch}{1.25}
	\setlength\arraycolsep{5pt}
	\begin{equation*}
		B = \begin{bmatrix}
					0.6\phantom{0} & 0.2\phantom{0} & 0.1\phantom{0} & 0.1\phantom{0}\\
					\cdot & 0.7\phantom{0} & 0.05\phantom{0} & 0.05\phantom{0}\\
					\cdot & \cdot & 0.6\phantom{0} & 0.25\phantom{0}\\
					\cdot & \cdot & \cdot & 0.6\phantom{0}
		\end{bmatrix}, \qquad
	\pi = (0.25, 0.25, 0.25, 0.25)^{\T},
	\end{equation*}
	\endgroup
	inspired by the simulation study in \cite{chen2020new:arxiv:v1}. We separately embed each graph adjacency matrix $A^{(j)}$ into dimension $k=4$ to obtain a varimax-transformed matrix estimate $\wh{Z}^{(j)} \in \R^{n \times k}$ using the function \textsf{varimax} in the \textsf{R} package \verb*|stats|. Given each graph embedding, we compute estimated class-conditional covariance matrices $\wh{\Sigma}^{(j)}_{\ell}$ for each cluster label $\ell \in \{1,2,3,4\}$ by fitting a one-component Gaussian density to the collection of vectors $(n\rho)^{1/2} \left\{\Pi_{Z} \cdot \wh{Z}^{(j)}_{i} - \breve{Z}_{i} \right\}$ satisfying $Z_{i} = e_{\ell}$. Our implementation uses the expectation--maximization algorithm found in the \textsf{R} package \verb*|mclust|. In this simulation setting, the signed permutation matrix $\Pi_{Z}$ can be computed using knowledge of the true latent factor matrix $Z$, though it is only needed to ensure proper ordering of coordinates and signs for comparison with the theoretical asymptotic covariance matrices. Here, \cref{table:example_fourblock_estimated_covariaces} reports the block-wise averages of the relative errors $\|\wh{\Sigma}^{(j)}_{\ell} - \Sigma_{\ell}\|_{\operatorname{F}}/\|\Sigma_{\ell}\|_{\operatorname{F}}$ computed in Frobenius norm, $1 \le \ell \le k$, together with sample standard errors, each scaled by a factor of one hundred for readability. As a further illustration, compare for example $\wh{\Sigma}^{(1)}_{1}$ with $\Sigma_{1}$ and $\wh{\Sigma}^{(1)}_{2}$ with $\Sigma_{2}$ below.

	\begingroup
	\renewcommand*{\arraystretch}{1.25}
	\setlength\arraycolsep{5pt}
	\begin{align*}
		\Sigma_{1} \overset{\cdot}{=}
		\begin{bmatrix} \phantom{-}14.7 & -5.84 & -1.71 & -1.71 \\ \cdot & \phantom{-}7.42 & \phantom{-}0.196 & \phantom{-}0.196 \\ \cdot & \cdot & \phantom{-}7.11 & -4.65 \\ \cdot & \cdot & \cdot & \phantom{-}7.11 \end{bmatrix}
		,&\hspace{1em}
		\wh{\Sigma}^{(1)}_{1} \overset{\cdot}{=}
		\begin{bmatrix} \phantom{-}14.4 & -6.00 & -1.75 & -1.69 \\ \cdot & \phantom{-}7.75 & \phantom{-}0.372 & \phantom{-}0.176 \\ \cdot & \cdot & \phantom{-}6.67 & -4.29 \\ \cdot & \cdot & \cdot & \phantom{-}6.72 \end{bmatrix}
		,\\
		\Sigma_{2} \overset{\cdot}{=}
		\begin{bmatrix} \phantom{-}10.2 & -5.25 & -1.07 & -1.07 \\ \cdot & \phantom{-}9.08 & \phantom{-}0.054 & \phantom{-}0.054 \\ \cdot & \cdot & \phantom{-}3.77 & -2.43 \\ \cdot & \cdot & \cdot & \phantom{-}3.77 \end{bmatrix}
		,&\hspace{1em}
		\wh{\Sigma}^{(1)}_{2} \overset{\cdot}{=}
		\begin{bmatrix} \phantom{-}11.4 & -5.92 & -1.04 & -1.43 \\\cdot & \phantom{-}9.41 & \phantom{-}0.215 & \phantom{-}0.079 \\ \cdot & \cdot & \phantom{-}3.32 & -2.10 \\ \cdot & \cdot & \cdot & \phantom{-}3.64 \end{bmatrix}.
	\end{align*}
	\endgroup
	
	Strictly speaking, for each stochastic blockmodel graph, \cref{thrm:main_SBM} does not simultaneously hold jointly for all order $\operatorname{O}(n)$ embedding vectors in each block $\ell$. In particular, there is dependence in the order $\operatorname{O}(n)$ vectors used to construct $\wh{\Sigma}^{(j)}_{\ell}$, and we do not claim that the estimation approach used in the current numerical example is optimal for estimating the asymptotic variances and covariances. Nevertheless, embedding the graphs and then aggregating the node embeddings follows conventional spectral graph clustering methodology, and doing so produces reasonably close agreement between empirical and theoretical covariances as shown in \cref{table:example_fourblock_estimated_covariaces}.

	\begin{table}[h]
	\def~{\hphantom{0}}
	\tbl{Estimation of asymptotic covariance matrices in \cref{sec:numerics_undirected}}{%
		\begin{tabular}{lcccc}
			\hline
			& Block 1 & Block 2 &  Block 3 & Block 4 \\[5pt]
			Rel.~Frob.~error (sample stand.~err.) $\times 100$ & 6.198 (0.2491) & 6.383 (0.2505) & 5.506 (0.2352) & 5.424 (0.2460)\\ \hline
	\end{tabular}}
	\label{table:example_fourblock_estimated_covariaces}
\end{table}

\subsection{Two-block directed stochastic blockmodel graphs}
\label{sec:numerics_directed}

\begin{figure}[t]%
		\centering
		\subfloat[\centering Vectors $\wh{Z}_{i}$ when $n=10^{3}$.]{{\includegraphics[width=4.5cm]{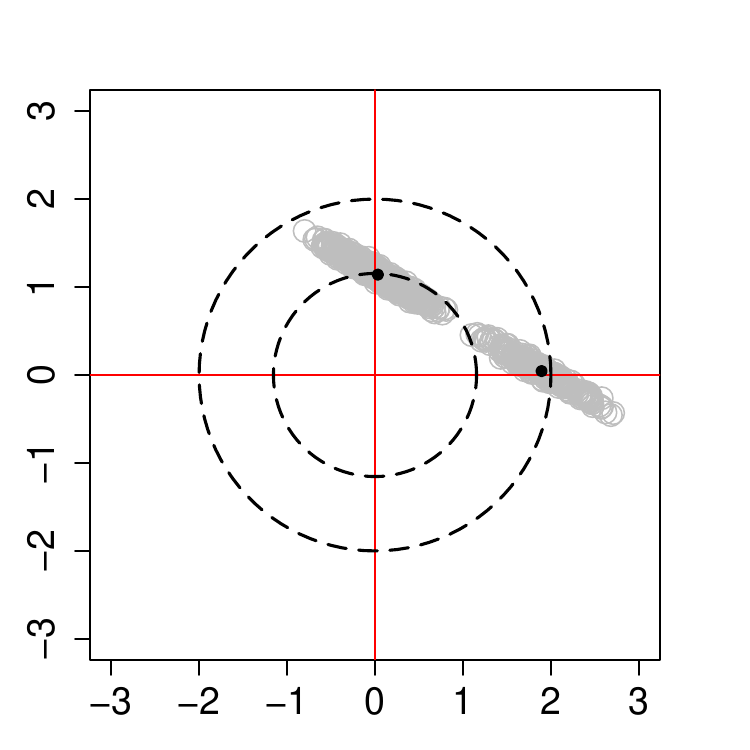} }}%
		\qquad
		\subfloat[\centering Vectors $\wh{Y}_{i}$ when $n=10^{3}$.]{{\includegraphics[width=4.5cm]{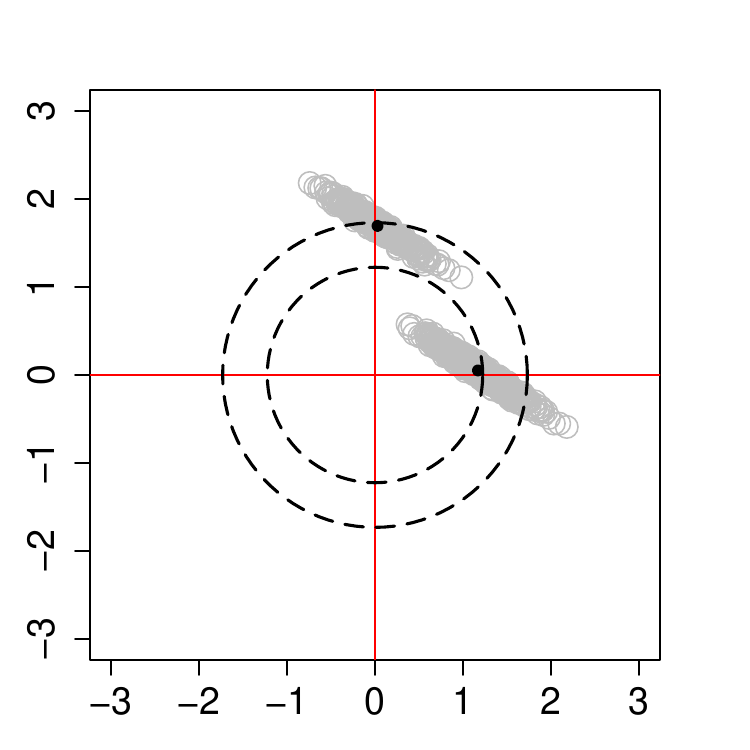} }} \\ \vspace{-1em}
		\subfloat[\centering Vectors $\wh{Z}_{i}$ when $n=10^{4}$.]{{\includegraphics[width=4.5cm]{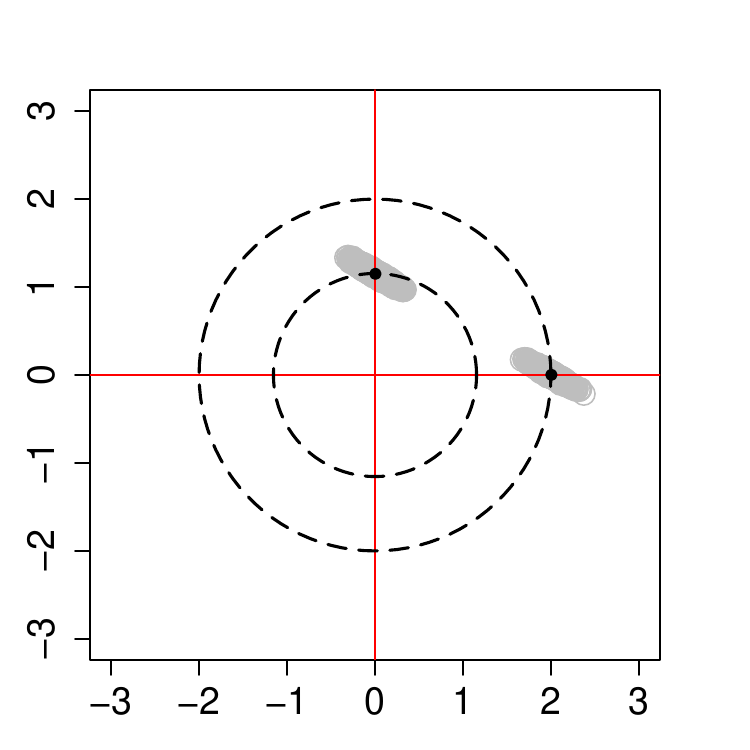} }}%
		\qquad
		\subfloat[\centering Vectors $\wh{Y}_{i}$ when $n=10^{4}$.]{{\includegraphics[width=4.5cm]{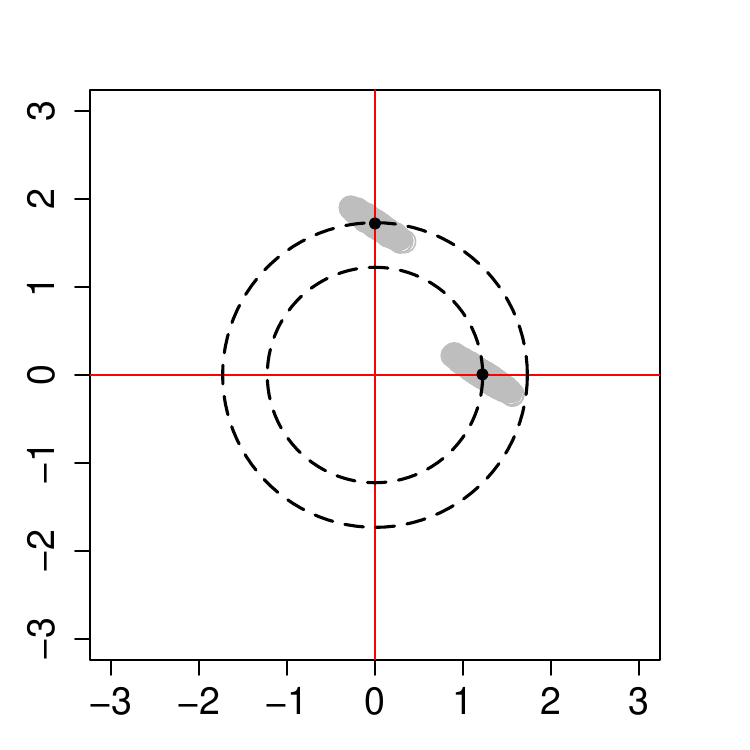} }}%
		\caption{Varimax-rotated left and right embeddings for two directed graphs. Embedding vectors are uncentered. Latent vectors are scaled in agreement with \cref{sec:main_DCSBM}. Dotted lines denote circles with radii given by the reciprocals of the square-roots of the block membership probabilities. Solid black circles denote the sample mean vectors computed for each block. Additional details are provided in \cref{sec:numerics_directed}.}
		\label{fig:example_directedSBM}
	\end{figure}

	Here, we simulate directed stochastic blockmodel graphs according to \cref{sec:main_DirSBM}, with left block membership vector $\pi_{Z} = (0.25, 0.75)^{\T}$, right block membership vector $\pi_{Y} = (0.667, 0.333)^{\T}$, sparsity term $\rho=1$, and block connectivity matrix defined entrywise as $B_{11} = 0.4, B_{12} = 0.6, B_{21} = 0.3, B_{22} = 0.7$. For two independent graphs of size $n=10^{3}$ and $n=10^{4}$, respectively, \cref{fig:example_directedSBM} depicts the varimax-transformed graph embeddings for the left and right estimated latent factor matrices. Crucially, the point cloud centroids lie along the standard coordinate axes, however each block exhibits elliptic radial covariance structure whose axes do not simultaneously align with the standard axes. In addition, each theoretical covariance matrix is reasonably well estimated in relative Frobenius norm by its sample counterpart.
	
\subsection{Degree-corrected stochastic blockmodels with affinity network structure}
\label{sec:numerics_degree-corrected}
	
\begin{figure}[t]%
		\centering
		\subfloat[\centering View of dimensions one and two]{{\includegraphics[width=4.5cm]{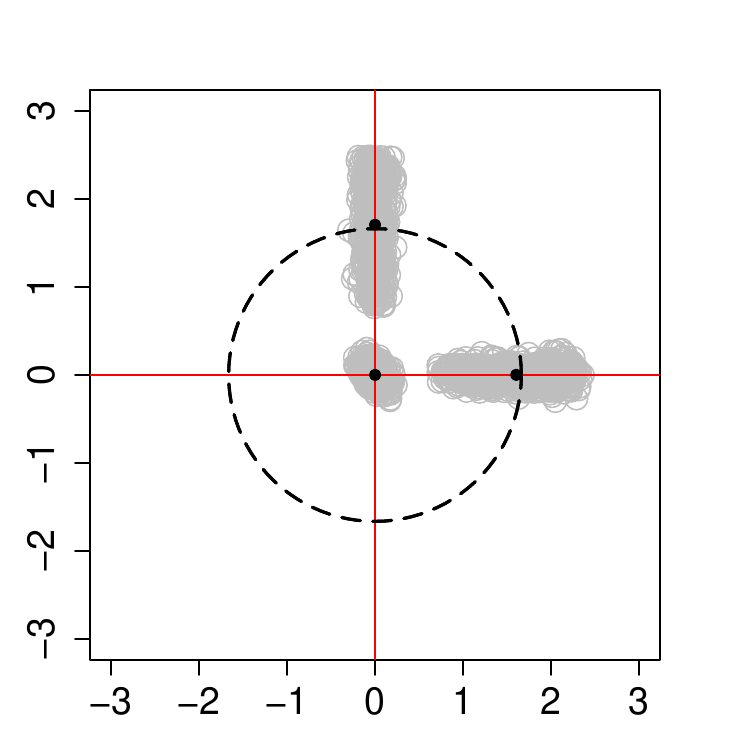} }}%
		\qquad
		\subfloat[\centering View of dimensions two and three]{{\includegraphics[width=4.5cm]{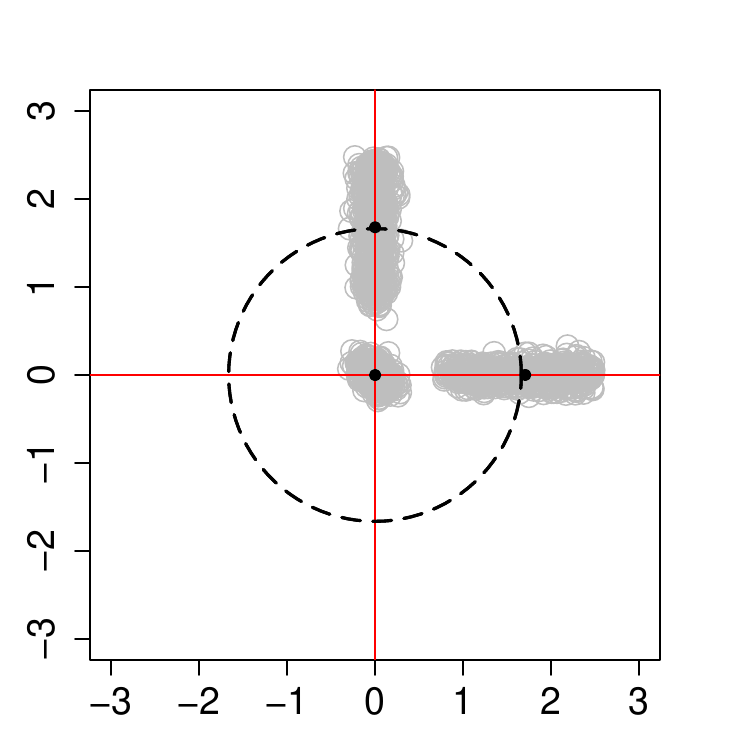} }}%
		\caption{Varimax-rotated degree-corrected stochastic blockmodel embedding for one graph. Dotted circles have radius equal to the magnitude of each block's non-zero theoretical coordinate centroid, while solid black circles denote empirical coordinate centroids. Additional details are provided in \cref{sec:numerics_degree-corrected}.}
		\label{fig:example_DCSBM}
	\end{figure}

	Here, we simulate a single dense three-block degree-corrected stochastic blockmodel graph per \cref{sec:main_DCSBM}, with two thousand nodes, equal membership probabilities, and block connectivity matrix defined entrywise as $B_{ii} = 0.2$ for all $1 \le i \le 3$ and $B_{ij} = 0.1$ for all $i \neq j$. The within-block edge probabilities are larger than the between-block edge probabilities which can be interpreted as producing affinity network structure. The node-specific degree parameters are drawn independently in the manner $\theta_{i} \sim \operatorname{Uniform}[0.25, 0.75]$ for all $1 \le i \le n$. For each panel in \cref{fig:example_DCSBM}, two of the three latent factor dimensions are isolated along the standard coordinate axes, with the one remaining dimension yielding coordinate values near the origin. This behavior generalizes to the same model having $k \ge 3$ blocks, illustrating the concept of ``simple structure'' \citep{thurstone1947multiple} as well as the concept of approximate sparsity in embedding space after performing dimensionality reduction.

\subsection{Asymptotic variance and covariance expressions}
\label{sec:numerics_asymp_covariance_examples}
	The asymptotic covariance expressions in \cref{thrm:main_SBM,thrm:main_DirSBM,thrm:main_DCSBM} are rather complicated in general. For simplicity and concreteness, consider undirected stochastic blockmodel graphs generated from the three-parameter model
	\begingroup
	\renewcommand*{\arraystretch}{1.25}
	\setlength\arraycolsep{5pt}
	\begin{equation*}
		B = \begin{bmatrix}
			a & b\\
			\cdot & a
		\end{bmatrix}, \quad 0 < b < a < 1; \qquad \pi=(\pi_{1},1-\pi_{1})^{\T}, \quad 0< \pi_{1} < 1.
	\end{equation*}
	\endgroup
	In this setting, it happens to be possible to analytically compute the block-specific asymptotic covariance matrices. In particular, for the first block, namely for the first community of nodes, a calculation reveals that
	\begingroup
	\renewcommand*{\arraystretch}{2.5}
	\setlength\arraycolsep{5pt}
	\begin{equation*}
		\Sigma_{1} =
		\begin{bmatrix}
			\frac{a^{3}(1-a)(1-\pi_{1}) + b^{3}(1-b)\pi_{1}}{(a^{2}-b^{2})^{2} \pi_{1}^{2}(1-\pi_{1})} & \frac{-a b \{a(1-a) + (a-b)(a+b-1)\pi_{1}\}}{(a^{2}-b^{2})^{2}\{\pi_{1}(1-\pi_{1})\}^{3/2}}\\
			\cdot & \frac{a^{2}b \pi_{1} + ab^{2}(1-a-\pi_{1})}{(a^{2}-b^{2})^{2}\pi_{1}(1-\pi_{1})^{2}}
		\end{bmatrix}.
	\end{equation*}
	\endgroup
	For the second block, a similar calculation reveals that
	\begingroup
	\renewcommand*{\arraystretch}{2.5}
	\setlength\arraycolsep{5pt}
	\begin{equation*}
		\Sigma_{2} = 
		\begin{bmatrix}
			\frac{ab^{2}\pi_{1} + a^{2}b(1-b-\pi_{1})}{(a^{2}-b^{2})^{2}\pi_{1}^{2}(1-\pi_{1})} & \frac{-ab\left\{b(1-b)+(a-b)(1-a-b)\pi_{1}\right\}}{(a^{2}-b^{2})^{2}\left\{\pi_{1}(1-\pi_{1})\right\}^{3/2}}\\
			\cdot & \frac{a^{3}(1-a)\pi_{1} + b^{3}(1-b)(1-\pi_{1})}{(a^{2}-b^{2})^{2}\pi_{1}(1-\pi_{1})^{2}}
		\end{bmatrix}.
	\end{equation*}
	\endgroup
	The determinants of the above covariance matrices, namely the generalized variances of the respective distributions, are given by
	\begin{equation*}
		\det(\Sigma_{1}) = \det(\Sigma_{2}) = \frac{a(1-a) \cdot b(1-b)}{(a^{2}-b^{2})^{2} \cdot \pi_{1}^{2}(1-\pi_{1})^{2}}.
	\end{equation*}
	The above expression succinctly reflects the near-degeneracy of the asymptotic covariance matrices in situations where the underlying stochastic blockmodel is itself nearly degenerate, namely when either $b \approx 0$, $b \approx a$, $a \approx 1$, $\pi_{1} \approx 0$, or $\pi_{1}\approx 1$.
	
	Evaluating the above expressions for several choices of model parameters yields the following examples of diverse covariance structure.
	\begin{align*}
		&a=\tfrac{3}{4}, b=\tfrac{1}{4}, \pi_{1}=\tfrac{1}{4}, &&\qquad \Sigma_{1} = \begin{bmatrix} 7 & -3^{1/2}\\ \cdot & \phantom{-}1\phantom{^{1/2}} \end{bmatrix}, &&\Sigma_{2} = \Sigma_{1};\\[1em]
		&a=\tfrac{4}{5}, b=\tfrac{3}{5}, \pi_{1}=\tfrac{1}{2}, &&\qquad \Sigma_{1} \overset{\cdot}{=} \begin{bmatrix} 9.63 & -9.79\\ \cdot & \phantom{-}10.77 \end{bmatrix}, &&\Sigma_{2} \overset{\cdot}{=} \begin{bmatrix} 10.77 & -9.79\\ \cdot & \phantom{-}9.63 \end{bmatrix};\\[1em]
		&a=\tfrac{8}{10}, b=\tfrac{3}{10}, \pi_{1}=\tfrac{4}{10}, &&\qquad \Sigma_{1} \overset{\cdot}{=} \begin{bmatrix} 2.37 & -1.21 \\ \cdot & \phantom{-}1.43 \end{bmatrix}, &&\Sigma_{2}  \overset{\cdot}{=} \begin{bmatrix} 2.97 & -1.28\\ \cdot & \phantom{-}1.20\end{bmatrix}.
	\end{align*}

\section{Extensions}
\label{sec:extensions}

	This paper considers stochastic blockmodel graphs having full-rank $k \times k$ connectivity matrix $B$ which is a widespread assumption in the literature. If instead $\operatorname{rank}(B) = r < k$, then modified versions of the results in \cref{sec:main} possibly still hold, though in such cases the $k$-dimensional singular subspaces of the latent factor matrices $Z$ and $Y$ are mapped into lower $r$-dimensional space for the population-level expectation term $Z B Y^{\T}$, resulting in quantitatively different varimax factor rotations and estimation properties. Broadly speaking, understanding such rank-degenerate blockmodeling is desirable because it permits flexibility in the true number of embedding dimensions $r$ and the true number of blocks or clusters $k$.
	
	This paper establishes the asymptotic normality of individual vectors representing node embeddings. Recent related work establishes that standardized averages of spectral estimates for the entries of $B$ also exhibit asymptotic normality but with bias in the presence of rank degeneracy that is exacerbated by network sparsity \citep{tang2022asymptotically}. Similar behavior is anticipated in the present paper when estimating the non-zero components of membership vectors via the standardized block-conditional averages of the vectors $\wh{Z}_{i}^{(j)}$; indeed, for the data in \cref{sec:numerics_undirected}, conducting a per-block Shapiro--Wilk test for normality yields $\operatorname{p}$-values larger than $0.2$.
	
	The findings in \cref{sec:main} possibly extend to non-blockmodel special cases of the more expressive generalized random dot product graph setting in \cite{rubindelanchy2022statistical} which permits more flexible latent variable modeling of networks beyond stochastic blockmodels. A challenge with such general, flexible models regarding varimax factor rotations, however, is the typical absence of even approximately sparse latent factor structure and the difficulty of establishing identifiability conditions for estimation.
	
	\cite{rohe2023vintage} develops theory and methods for more general data-generating models than stochastic blockmodels. Additional steps such as data centering, recentering, and row-normalized varimax are discussed therein. In future work, it would be interesting to quantify the asymptotic distributional properties of varimax-based estimators in such settings, which foreseeably would introduce additional complexities and technical challenges.

\section{Discussion}
\label{sec:discussion}

	Stochastic blockmodels have attracted significant research interest in recent years at the time of writing this paper. Despite their seemingly straightforward structure, the statistics community is only now developing a more complete understanding of spectral-based inference for these baseline statistical network models \citep{rohe2011spectral,lei2015consistency,lei2016goodness,cape2019signal,tang2022asymptotically,athreya2022eigenvalues}. While stochastic blockmodels are of moderate interest in and of themselves, studying them ultimately provides points of departure for investigating more complex relational data and for moving beyond node clustering or the so-called community detection problem.
	
	Readers familiar with the statistical network analysis literature on spectral embedding will notice that the results in this paper somewhat resemble those found in several existing works, e.g.,~\cite{athreya2016limit,tang2018limit,rubindelanchy2022statistical}. These related works and others primarily establish the uniform asymptotic recovery of node embeddings and asymptotic normality for sequences of vectors acted on by unknown linear transformations, typically indefinite or orthogonal transformations, which admit no general convergence properties or obvious interpretability. Consequently, such results are more useful for informing clustering methodology than for direct estimation in practice. In contrast, our results incorporate varimax factor rotations with graph embeddings to enable more straightforward estimation modulo only the set of signed permutation matrices $\signorthGroup(\cdot)$. Advantageously, the collection $\signorthGroup(\cdot)$ reflects the minimal possible necessary statistical non-uniqueness in general, since singular vectors in left-right pairs are uniquely computed only up to a choice of sign and since blockmodels are unique only up to globally permuting the latent cluster labels or the block coordinates.
	
	In the undirected and degree-corrected stochastic blockmodel numerical examples, nodes within the same graph share a common expected degree, hence preventing trivial degree-based node clustering. In the directed stochastic blockmodel numerical example, separately applying varimax to the estimated left and right singular vectors is shown to successfully recover left-right directional membership heterogeneity. Our simulations make use of the true latent block memberships which are asymptotically perfectly recoverable via $k$-means clustering in the present semi-sparse regime $n \rho = \omega(\log^{c} n)$. Moreover, the embedding dimension can be asymptotically perfectly recovered via an elbow in the scree plot, namely via a demonstrable spectral gap.
	
	This paper considers computing a truncated singular value decomposition of a given matrix $A$ to obtain estimated latent factor matrices $\wh{Z}$ and $\wh{Y}$ for $Z$ and $Y$ appearing in $\Ex(A \mid Z, Y)$, per \cref{eq:semiparFM_expectation}. As such, the methods and results in this paper correspond to the so-called problem of matrix denoising. Matrix denoising is distinct from the problem of covariance estimation in factor models via principal components analysis. In particular, this paper allows for heteroskedasticity in the underlying noise matrix, whereas it is well-known that heteroskedastic noise can pose significant challenges for classical principal components analysis and factor analysis. Further discussion and results to this effect can be found in, for example, \cite{zhang2022heteroskedastic,yan2021inference}.
	
	Unanswered questions abound. For example, the statistical properties of other common factor rotation methods, such as promax, quartimax, and equimax, are not fully developed, even for semiparametric factor models. The present paper does not consider the optional row-normalization step in varimax which remains popular with practitioners and whose theoretical properties are incompletely understood. For networks specifically, this paper considers low-dimensional embeddings derived from adjacency matrices, while the widespread success and popularity of graph Laplacians makes them natural candidates to consider together with varimax rotations. In this direction, it has been shown that Laplacian-based analyses perform well empirically, even for high-dimensional, relatively sparse data matrices, e.g.,~see \cite{rohe2023vintage,chen2020new:arxiv:v1}. That being said, their theoretical properties can be quite challenging to analyze, e.g.,~see \cite{tang2018limit}. These are but a few of the many open problems in contemporary multivariate analysis.

\section*{Acknowledgement}
\label{sec:acknowledgement}

	The author thanks the Editor, the Associate Editor, and two referees for their valuable feedback and suggestions. The author also thanks Karl Rohe and Minh Tang for stimulating discussions. This research is supported in part by the National Science Foundation under grants DMS-1902755 and SES-1951005. The author gratefully acknowledges support from the University of Wisconsin--Madison, Office of the Vice Chancellor for Research and Graduate Education with funding from the Wisconsin Alumni Research Foundation.


\bibliographystyle{biometrika}
\bibliography{bib}



\newpage
\section{Supplementary Material} 


\subsection{Asymptotic normality in varimax-transformed undirected graph embeddings}
\label{sec:proof_vmx_undirected_sbm}

\begin{proof}[of Theorem~1]
	
	This proof proceeds via matrix eigendecompositions, or equivalently, spectral decompositions. This choice is primarily for convenience; we emphasize that the analysis also holds for matrix singular value decompositions, with only minor notational modifications, since the matrices $A$ and $B$ of primary interest are real-valued and symmetric.
	
	Given the full-rank symmetric matrix of probabilities $B \in (0, 1)^{k \times k}$, write its spectral decomposition as
	\begin{equation*}
		B = U_{B} \Lambda_{B} U_{B}^{\T} = U_{B} \cdot \left(|\Lambda_{B}|^{1/2} \cdot J \cdot |\Lambda_{B}|^{1/2}\right) \cdot U_{B}^{\T},
	\end{equation*}
	where the $k \times k$ orthogonal matrix $U_{B}$ consists of columns that are unit-norm eigenvectors of $B$, the $k \times k$ diagonal matrix $\Lambda_{B}$ contains the eigenvalues of $B$ arranged in decreasing order, $|\cdot|$ denotes the entrywise absolute value operation, and $J$ is the $k \times k$ diagonal matrix $J = \diag(1_{p}, -1_{q})$ when $B$ has $p$ positive eigenvalues and $q$ negative eigenvalues, hence $p+q=k$. The matrix $U_{B}$ is not uniquely defined, modulo column-wise signflips and up to subspace basis non-uniqueness if $B$ has repeated eigenvalues, but this non-uniqueness poses no appreciable problem throughout the proof due to appropriate invariances, the consideration of certain matrix products, and the bounding of subspace distances.
	
	By hypothesis, the sparse latent factor matrix $Z \in \{0,1\}^{n \times k}$ consists of independent, identically distributed $k$-dimensional random latent factor vectors $Z_{i} \sim \mathcal{G} = \operatorname{Multinomial}(1; \pi)$ whence $Z_{i} \in \{e_{1},\dots,e_{k}\}$ and where $\pi \in \R^{k}$ is a fixed probability vector in the interior of the unit simplex. Now, define the matrix of coefficients $T_{B} = |\Lambda_{B}|^{1/2} \cdot U_{B}^{\T}$, and for each index $i$ define the $k$-dimensional random vector $\xi_{i} = T_{B} \cdot Z_{i}$. Write $X = [\xi_{1},\dots,\xi_{n}]^{\T} \in \R^{n \times k}$ as the matrix whose rows are the vectors $\xi_{i}$, hence $Z B Z^{\T} = X J X^{\T}$ almost surely. At the population level, $\Ex(A \mid Z) = \Ex(A \mid X)$.
	
	For semi-sparse undirected stochastic blockmodel graphs satisfying $n \rho = \omega(\log^{c} n)$ with positive-definite connectivity matrix $B$, it holds as in \cite[Theorem~3]{cape2019signal} that for any fixed index $i$, the sequence of random vectors of the form
	\begin{equation*}
		n\rho^{1/2} \left\{W_{Z}^{\T} \cdot \left(W_{U}\wh{U}_{i} - U_{i}\right) \mid Z_{i}=e_{\ell}\right\}
	\end{equation*}
	converges in distribution to a mean zero multivariate Gaussian random vector as $n \rightarrow \infty$, where (i) the explicit asymptotic covariance matrix depends on the block membership of the $i$-th node, indicated by $e_{\ell}$, and (ii) the orthogonal matrices $W_{Z}$ and $W_{U}$, implicitly indexed by $n$, form sequences denoted by $(W_{Z}), (W_{U}) \subset \orthGroup(k)$. Here, the column vector $\wh{U}_{i}$ is the $i$-th row of the matrix $\wh{U}$ consisting of the top-$k$ unit norm eigenvectors of $A$. Similarly, the column vector $U_{i}$ is the $i$-th row of the matrix $U$ consisting of the top-$k$ unit norm eigenvectors of $\Ex(A \mid Z)$. This existing result does not involve varimax factor rotations.
	
	In this proof, we shall first generalize the above asymptotic normality result to hold for all full-rank stochastic blockmodels, not only those having positive-definite connectivity matrices $B$. In other words, we shall first extend the existing result \cite[Theorem 3]{cape2019signal}; the minor distinction between fixed versus generated block memberships is immaterial here. Our extension will be achieved by revisiting and rewriting matrix factorizations involving $U$ and $\Lambda^{-1}$ as equivalent, properly oriented matrix factorizations involving $X$ and $J$. To this end, it will be convenient to work with the reparametrization and minor change of notation established at the start of this proof.
	
	Let $P$ denote the matrix $\Ex(A \mid Z)$, equivalently $\Ex(A \mid X)$, which in terms of $X$ is written as $P = \rho X J X^{\T}$ and which has rank $k$ almost surely. Denote the skinny spectral decomposition of $P$ by $U \Lambda U^{\T}$, whence it follows that $P^{\dagger} = U \Lambda^{-1} U^{\T}$ is the necessarily unique Moore--Penrose inverse of $P$. Consequently, by standard properties of the Moore--Penrose inverse,
	\begin{equation*}
		U \Lambda^{-1} U^{\T} = P^{\dagger} = (\rho X J X^{\T})^{\dagger} = \rho^{-1}(X^{\T})^{\dagger}J^{\dagger}X^{\dagger} =  \rho^{-1}X(X^{\T}X)^{-1} J (X^{\T}X)^{-1}X^{\T}.
	\end{equation*}
	Squaring the matrix $U \Lambda^{-1} U^{\T}$ yields
	\begin{equation*}
		U \Lambda^{-2} U^{\T} = \rho^{-2}X(X^{\T}X)^{-1}J(X^{\T}X)^{-1}J(X^{\T}X)^{-1}X^{\T},
	\end{equation*}
	which is a positive-semidefinite matrix, hence there exists a $k \times k$ orthogonal matrix $W_{X}$ such that
	\begin{equation*}
		U \Lambda^{-1} W_{X} = \rho^{-1}X(X^{\T}X)^{-1}J(X^{\T}X)^{-1/2}.
	\end{equation*}
	Let $E = A - P$, hence it follows that
	\begin{equation*}
		E U \Lambda^{-1} W_{U} = E \left( U \Lambda^{-1} W_{X} \right) W_{X}^{\T} W_{U} = \rho^{-1}EX(X^{\T}X)^{-1}J(X^{\T}X)^{-1/2}W_{X}^{\T}W_{U}.
	\end{equation*}
	By taking the above display equation and proceeding via the same perturbation analysis approach as in \cite[Supplement]{cape2019signal}, we obtain the key matrix identity
	\begin{equation*}
		\left(\wh{U}W_{U}^{\T} - U\right)W_{X} = \rho^{-1}EX(X^{\T}X)^{-1}J(X^{\T}X)^{-1/2} + E^{\prime}W_{U}^{\T}W_{X},
	\end{equation*}
	where $E^{\prime} \in \R^{n \times k}$ denotes a residual matrix of smaller order in maximum Euclidean row norm. Namely, the $i$-th row vector in the above display equation, when written as a column vector and scaled by $n \rho^{1/2}$, satisfies
	\begin{equation}
		\label{eq:prf_biom_vector_extension}
		n\rho^{1/2} W_{X}^{\T}\left(W_{U}\wh{U}_{i} - U_{i}\right) = (n^{-1}X^{\T}X)^{-1/2} J (n^{-1}X^{\T}X)^{-1} \left\{(n\rho)^{-1/2} (EX)_{i} \right\} + 
		\ohp(1),
	\end{equation}
	where the notation $\chi_{1} = \chi_{2}+ \ohp(1)$ conveys $\|\chi_{1} - \chi_{2}\|_{\ell_{2}} = \ohp(1)$ as $n \rightarrow \infty$. Hence, by an application of the multivariate central limit theorem and Slutsky's theorem, the sequence of random vectors of the form
	\begin{equation}
		\label{eq:prf_biometrika_clt_extension}
		n\rho^{1/2} \left\{W_{X}^{\T} \cdot \left(W_{U}\wh{U}_{i} - U_{i}\right) \mid Z_{i}=e_{\ell}\right\}
	\end{equation}
	converges in distribution to a multivariate normal random vector with mean zero and covariance matrix $\Delta_{X}^{-1/2} J \Delta_{X}^{-1} \cdot \Ex_{\xi}\{ g_{\rho_{\infty}}(\xi, T_{B} e_{\ell}) \} \cdot \Delta_{X}^{-1} J \Delta_{X}^{-1/2}$ as $n \rightarrow \infty$. Here, $\Delta_{X} = \Ex(\xi\xi^{\T}) \in \R^{k \times k}$ where $\xi \sim T_{B}\mathcal{G}$, while $g_{\rho_{\infty}}(\xi, x) = (x^{\T}J\xi) \cdot (1-\rho_{\infty}x^{\T}J\xi) \cdot \xi\xi^{\T}$. Recall that by assumption, $n \rho = \omega(\log^{c}n)$ with $\rho = \rho_{n} \in (0,1]$ and $\rho_{n} \rightarrow \rho_{\infty} \in \{0,1\}$. In words, the asymptotic covariance matrix differs between the dense graph regime $\rho = 1$ and the semi-sparse graph regime $\rho \rightarrow 0$, exhibiting larger asymptotic variability in the latter case.
	
	
	The next step of the proof is to synthesize the above perturbation analysis together with the effect of varimax factor rotations. Our approach will make use of existing strong consistency results, namely $\ell_{\ttinf}$ uniform row-wise error bounds. We shall also need to demonstrate suitable orthogonal alignment of certain matrix terms appearing in perturbation expressions.
	
	Towards this end, first observe that $U W_{X} = (U \Lambda U^{\T}) \cdot (U \Lambda^{-1} W_{X}) = X(X^{\T}X)^{-1/2}$, hence for each fixed index $i$ it holds almost surely that
	\begin{equation*}
		W_{X}^{\T} \cdot n^{1/2}U_{i}
		= (n^{-1}X^{\T}X)^{-1/2}X_{i}\\
		= (n^{-1}X^{\T}X)^{-1/2} \cdot T_{B} \cdot Z_{i}.
	\end{equation*}
	
	At the same time, $W_{U}$ denotes the closest $k$-dimensional orthogonal matrix approximation to $U^{\T}\wh{U}$ which satisfies the operator norm bound $\|U^{\T}\wh{U}  - W_{U}\|_{\op} = \Ohp\left\{(n\rho)^{-1}\right\}$ per \cite[Eq.~(S2) in Supplement]{cape2019signal}. Furthermore, applying \cite[Theorem~1]{cape2019signal} yields the maximum Euclidean row norm bound $\|\wh{U} - UW_{U}\|_{\ttinf} = \ohp\left(n^{-1/2}\right)$. Consequently, $n^{1/2}\|\wh{U}\|_{\ttinf} = \Ohp(1)$ since $n^{1/2}\|U\|_{\ttinf}$ is bounded with overwhelming probability. Thus, since $\wh{Z}_{i} = n^{1/2}R_{\wh{U}}^{\T}\wh{U}_{i}$ and $X = Z T_{B}^{\T}$, it follows that
	\begin{align*}
		& W_{X}^{\T}W_{U} \cdot n^{1/2}\wh{U}_{i}\\
		&\hspace{1em}= W_{X}^{\T}\left(U^{\T}\wh{U}\right) \cdot n^{1/2}\wh{U}_{i} + \ohp\left\{(n\rho)^{-1/2}\right\}\\
		&\hspace{1em}= \left(X^{\T}X\right)^{-1/2}X^{\T} \cdot \wh{U}\left(R_{\wh{U}} R_{\wh{U}}^{\T}\right) \cdot n^{1/2}\wh{U}_{i} + \ohp\left\{(n\rho)^{-1/2}\right\}\\
		&\hspace{1em}= \left(n^{-1}X^{\T}X\right)^{-1/2} \cdot T_{B} \cdot \left(n^{-1}Z^{\T}\wh{Z}\right) \cdot \wh{Z}_{i} + \ohp\left\{(n\rho)^{-1/2}\right\}.
	\end{align*}
	Given the asymptotic normality of the eigenvector components in \cref{eq:prf_biometrika_clt_extension}, an application of Slutsky's theorem together with several steps of algebra yields that
	\begin{equation*}
		(n\rho)^{1/2}\left\{\left(n^{-1}Z^{\T}\wh{Z}\right) \cdot \wh{Z}_{i} - Z_{i} \mid Z_{i} = e_{\ell} \right\}
	\end{equation*}
	is asymptotically normal with mean zero. In particular, necessarily $(n^{-1}Z^{\T}\wh{Z}) \cdot \wh{Z}_{i} \rightarrow Z_{i}$ in probability given $Z_{i} = e_{\ell}$, since otherwise we obtain a contradiction to the established convergence in distribution to a multivariate normal random vector.
	
	For any choice of $\Pi_{Z} \in \signorthGroup(k)$, adding and subtracting $Z \diag(\pi)^{-1/2} \Pi_{Z}$ below gives
	\begin{equation}
		\label{eq:prf_Zhat_Z_product_undirected}
		\left(n^{-1}Z^{\T}\wh{Z}\right)
		= \left(n^{-1}Z^{\T}Z\right)\diag(\pi)^{-1/2}\Pi_{Z} + n^{-1}Z^{\T}\left(\wh{Z} - Z\diag(\pi)^{-1/2}\Pi_{Z}\right).
	\end{equation}
	The law of large numbers guarantees that $\left(n^{-1}Z^{\T}Z\right) \rightarrow \diag(\pi)$ almost surely as $n \rightarrow \infty$. Furthermore, \cite[Corollary~C.1]{rohe2023vintage} establishes that there exists a sequence of signed permutation matrices $\left(\Pi_{Z}\right) \subset \signorthGroup(k)$ such that
	\begin{equation*}
		\max_{1 \le i \le n} \left\|\left(\wh{Z} - Z\diag(\pi)^{-1/2}\Pi_{Z}\right)_{i}\right\|_{\ell_{2}} = \ohp(1).
	\end{equation*}
	For this same sequence of signed permutation matrices, properties of the $\ell_{\ttinf}$ norm yield
	\begin{equation*}
		\left\|n^{-1}Z^{\T}\left(\wh{Z} - Z\diag(\pi)^{-1/2}\Pi_{Z}\right)\right\|_{\ttinf} \le \left\|n^{-1}Z^{\T}\right\|_{\inftyinfty} \cdot \left\|\wh{Z} - Z\diag(\pi)^{-1/2}\Pi_{Z}\right\|_{\ttinf},
	\end{equation*}
	hence,
	\begin{equation*}
		\left\|n^{-1}Z^{\T}\left(\wh{Z} - Z\diag(\pi)^{-1/2}\Pi_{Z}\right)\right\|_{\ttinf} = \ohp(1),
	\end{equation*}
	where $\|\cdot\|_{\ttinf}$ again denotes the maximum Euclidean row norm of a matrix, where $\|\cdot\|_{\inftyinfty}$ denotes the maximum absolute row sum of a matrix, and where we have used the fact that $\left\|n^{-1}Z^{\T}\right\|_{\inftyinfty} \le 1$ almost surely. Moreover, $\Pi_{Z} \wh{Z}_{i} \rightarrow \diag(\pi)^{-1/2} Z_{i}$ in probability, recalling that $\diag(\pi)^{-1/2} Z_{i} = \breve{Z}_{i}$ with $Z_{i}=e_{\ell}$. Taken together and by writing $Z_{i} = \diag(\pi)^{1/2}\diag(\pi)^{-1/2} Z_{i}$, these deductions imply that
	\begin{equation*}
		(n\rho)^{1/2}\left\{\Pi_{Z} \cdot \wh{Z}_{i} - \breve{Z}_{i} \mid Z_{i} = e_{\ell} \right\}
	\end{equation*}
	is asymptotically multivariate normal as $n \rightarrow \infty$. By tracing through the proof and again appealing to Slutsky's theorem, the asymptotic covariance matrix is given by
	\begin{equation*}
		\left\{\diag(\pi)^{-1/2} \cdot T_{B}^{-1} \cdot J \cdot \Delta_{X}^{-1}\right\} \cdot \Ex_{\xi}\{ g_{\rho_{\infty}}(\xi, T_{B} e_{\ell}) \} \cdot \left\{\Delta_{X}^{-1} \cdot J \cdot (T_{B}^{-1})^{\T} \cdot \diag(\pi)^{-1/2}\right\}.
	\end{equation*}
	To reiterate, here $\pi$ is the block membership probability vector for the latent factor vectors, here $T_{B} = |\Lambda_{B}|^{1/2} U_{B}^{\T}$ is defined such that $B = T_{B}^{\T} J T_{B}$ where $J = \diag(1_{p},-1_{q}) \in \R^{k \times k}$, here $\Delta_{X} = \Ex(\xi\xi^{\T})$ where $\xi$ has distribution given by $T_{B} \mathcal{G}$, and $g_{\rho_{\infty}}(\xi, x) = (x^{\T}J\xi) \cdot (1-\rho_{\infty}x^{\T}J\xi) \cdot \xi\xi^{\T}$ where $\rho_{\infty} \in \{0,1\}$. This concludes the proof.
\end{proof}

\subsection{Asymptotic normality in varimax-transformed directed graph embeddings}
\label{sec:proof_vmx_directed_sbm}
\begin{proof}[of Theorem~2]
	
	By hypothesis, consider directed stochastic blockmodels having edge probabilities given by the matrix $\Ex(A \mid Z, Y) = \rho Z B Y^{\T}$, consisting of independent, identically distributed left latent factor vectors $Z_{i} \sim \mathcal{G} = \operatorname{Multinomial}(1; \pi_{Z})$ and independent, identically distributed right latent factor vectors $Y_{i} \sim \mathcal{H} = \operatorname{Multinomial}(1; \pi_{Y})$ and where $Z_{i}$ is independent of $Y_{i}$ for all $1 \le i \le n$. Write the singular value decomposition of the full-rank $k \times k$ matrix $B$ as $B = U_{B} S_{B} V_{B}^{\T}$. Define the random vectors $\xi_{i} = T_{\Xi} \cdot Z_{i} = S_{B}^{1/2}U_{B}^{\T} \cdot Z_{i}$ and $\upsilon_{i} = T_{\Upsilon} \cdot Y_{i} = S_{B}^{1/2}V_{B}^{\T} \cdot Y_{i}$ to form the rows of $\Xi$ and $\Upsilon$, whence $Z B Y^{\T} = \Xi \Upsilon^{\T}$ almost surely. Consequently, $\Ex(A \mid Z,Y) = \Ex(A \mid \Xi, \Upsilon)$ which for convenience shall be denoted simply by $P$.
	
	By extending the proof analysis found in \cref{sec:proof_vmx_undirected_sbm} and in the supplementary material of \cite{cape2019signal} to the case of directed graphs and their truncated singular value decompositions, we obtain the key expression 
	\begin{equation}
		\label{eq:prf_biom_vector_extension_directed}
		n\rho^{1/2} \cdot W_{\Upsilon}^{\T} \left( W_{V}\wh{U}_{i} - W_{V}W_{U}^{\T} U_{i}\right) = (n^{-1}\Xi^{\T}\Xi)^{-1/2} (n^{-1}\Upsilon^{\T}\Upsilon)^{-1}\left\{(n\rho)^{-1/2} (E\Upsilon)_{i}\right\} + \ohp(1)
	\end{equation}
	which generalizes \cref{eq:prf_biom_vector_extension}. \cref{eq:prf_biom_vector_extension_directed} can alternatively be established with the help of Theorem~9 in \cite{yan2021inference}. In the above display equation, $W_{U}$ minimizes $\|\wh{U}-UW\|_{\Frob}$ over all $k \times k$ orthogonal matrices $W$. Separately, $W_{V}$ minimizes $\|\wh{V}-VW\|_{\Frob}$ over all $k \times k$ orthogonal matrices $W$. Here, $W_{\Upsilon}$ denotes the orthogonal matrix satisfying
	\begin{equation*}
		V S^{-1} W_{\Upsilon} = \rho^{-1}(\Upsilon^{\T}\Upsilon)^{-1}(\Xi^{\T}\Xi)^{-1/2},
	\end{equation*}
	hence $U W_{\Upsilon} = (USV^{\T}) \cdot (V S^{-1} W_{\Upsilon}) = \Xi (\Xi^{\T}\Xi)^{-1/2}$. By analogy, for the orthogonal matrix $W_{\Xi}$, it holds that $U S^{-1} W_{\Xi} = \rho^{-1}(\Xi^{\T}\Xi)^{-1}(\Upsilon^{\T}\Upsilon)^{-1/2}$ and so $V W_{\Xi} = \Upsilon(\Upsilon^{\T}\Upsilon)^{-1/2}$.
	
	A routine computation shows $\|U^{\T}\wh{U} - W_{U}\|_{\op} \le \|\wh{U}\wh{U}^{\T} - UU^{\T}\|_{\op}^{2} = \Ohp\left\{(n\rho)^{-1}\right\}$ and similarly $\|V^{\T}\wh{V} - W_{V}\|_{\op} = \Ohp\left\{(n\rho)^{-1}\right\}$, where the stated bounds hold by an application of Wedin's theorem, see for example \citep[Theorem~2.9]{chen2021spectral}.
	
	Next, consider the expansion
	\begin{align*}
		& W_{\Upsilon}^{\T}W_{V}\\
		&\hspace{1em}= W_{\Upsilon}^{\T}V^{\T}\wh{V} + \ohp\left\{(n\rho)^{-1/2}\right\}\\
		&\hspace{1em}= W_{\Upsilon}^{\T} W_{\Xi} W_{\Xi}^{\T} V^{\T}\wh{V} + \ohp\left\{(n\rho)^{-1/2}\right\}\\
		&\hspace{1em}= W_{\Upsilon}^{\T} W_{\Xi} (\Upsilon^{\T}\Upsilon)^{-1/2}\Upsilon^{\T}\wh{V} + \ohp\left\{(n\rho)^{-1/2}\right\}\\
		&\hspace{1em}= W_{\Upsilon}^{\T} W_{\Xi} \cdot  (n^{-1}\Upsilon^{\T}\Upsilon)^{-1/2} \cdot T_{\Upsilon} \cdot n^{-1} \left(Y^{\T} \cdot n^{1/2}\wh{V}R_{\wh{V}}\right)R_{\wh{V}}^{\T} + \ohp\left\{(n\rho)^{-1/2}\right\}\\
		&\hspace{1em}= W_{\Upsilon}^{\T} W_{\Xi} \cdot  (n^{-1}\Upsilon^{\T}\Upsilon)^{-1/2} \cdot T_{\Upsilon} \cdot \left(n^{-1} Y^{\T} \wh{Y}\right)R_{\wh{V}}^{\T} + \ohp\left\{(n\rho)^{-1/2}\right\}.
	\end{align*}
	Crucially, on the right-hand side, the leftmost matrix product can be expanded as
	\begin{align*}
		&W_{\Upsilon}^{\T}W_{\Xi}\\
		&\hspace{1em}= W_{\Upsilon}^{\T} U^{\T}U V^{\T}V W_{\Xi}\\
		&\hspace{1em}= (\Xi^{\T}\Xi)^{-1/2}\Xi^{\T} \cdot \left(U V^{\T}\right) \cdot \Upsilon(\Upsilon^{\T}\Upsilon)^{-1/2}\\
		&\hspace{1em}= (\Xi^{\T}\Xi)^{-1/2} \cdot (T_{\Xi} Z^{\T}) \cdot \left(U V^{\T}\right) \cdot (Y T_{\Upsilon}^{\T}) \cdot (\Upsilon^{\T}\Upsilon)^{-1/2} \\
		&\hspace{1em}= (\Xi^{\T}\Xi)^{-1/2} \cdot (T_{\Xi} Z^{\T}) \cdot \left(n^{-1} Z (Z^{\T}Z)^{-1/2} \wt{R}_{U}^{\T} \wt{R}_{V} (Y^{\T}Y)^{-1/2} Y^{\T}\right) \cdot (Y T_{\Upsilon}^{\T}) \cdot (\Upsilon^{\T}\Upsilon)^{-1/2},
	\end{align*}
	where the orthogonal matrix $\wt{R}_{U}^{\T}$ consists of the left singular vectors of the matrix $(n^{-1}Z^{\T}Z)^{1/2}B(n^{-1}Y^{\T}Y)^{1/2}$ and $\wt{R}_{V}$ consists of the corresponding right singular vectors.
	
	The matrix product $(n^{-1}Y^{\T}Y)^{-1/2}(n^{-1}Y^{\T}Y) T_{\Upsilon}^{\T}(n^{-1}\Upsilon^{\T}\Upsilon)^{-1/2}(n^{-1}\Upsilon^{\T}\Upsilon)^{-1/2}T_{\Upsilon} = (n^{-1}Y^{\T}Y)^{-1/2}$ converges to $\diag(\pi_{Y})^{-1/2}$ almost surely as $n \rightarrow \infty$ by the law of large numbers. Moreover, by the same proof strategy used to analyze \cref{eq:prf_Zhat_Z_product_undirected}, there exists a sequence of signed permutation matrices $(\Pi_{Y}) \subset \signorthGroup(k)$ indexed by $n$ such that $(n^{-1} Y^{\T} \wh{Y}) \Pi_{Y}^{\T}$ converges to $\diag(\pi_{Y})^{1/2}$. By putting the pieces together and with a slight abuse of notation,
	\begin{align*}
		& W_{\Upsilon}^{\T} W_{V} \cdot (n^{1/2} \wh{U}_{i})\\
		&\hspace{2em}= \left\{\Delta_{\Xi}^{-1/2} T_{\Xi} \diag(\pi_{Z})^{1/2}\right\}
		\cdot \left( \wt{R}_{U}^{\T} \wt{R}_{V} \Pi_{Y} R_{\wh{V}}^{\T} R_{\wh{U}}\right)
		\cdot \left(n^{1/2} R_{\wh{U}}^{\T} \wh{U}_{i}\right) + \ohp\left\{(n\rho)^{-1/2}\right\}.
	\end{align*}
	Above, we are using the fact that $\|\wh{U} - UW_{U}\|_{\ttinf} = \ohp\left(n^{-1/2}\right)$ and consequently $n^{1/2}\|\wh{U}\|_{\ttinf} = \Ohp(1)$ since $n^{1/2}\|U\|_{\ttinf}$ is bounded with overwhelming probability, as in the setting of undirected stochastic blockmodels. Similarly, $\|\wh{V} - VW_{V}\|_{\ttinf} = \ohp\left(n^{-1/2}\right)$ and $n^{1/2}\|\wh{V}\|_{\ttinf} = \Ohp(1)$.
	
	Next, observe that for any choice of sequence $(\Pi_{Z}) \subset \signorthGroup(k)$, for each value of $n$,
	\begin{equation*}
		\left(\wt{R}_{U}^{\T} \wt{R}_{V} \Pi_{Y} R_{\wh{V}}^{\T} R_{\wh{U}}\right) \Pi_{Z}^{\T}
		= \Pi_{Z} R_{\wh{U}}^{\T} \cdot \left(R_{\wh{U}}\Pi_{Z}^{\T}\wt{R}_{U}^{\T}\right) \cdot \left(R_{\wh{V}}\Pi_{Y}^{\T}\wt{R}_{V}^{\T}\right)^{\T} \cdot R_{\wh{U}} \Pi_{Z}^{\T},
	\end{equation*}
	hence choosing the sequences $(\Pi_{Z})$ and $(\Pi_{Y})$ per \cite{rohe2023vintage} yields that the matrix product in the preceding display equation converges in probability to the identity matrix in the large-$n$ limit.
	
	In a similar fashion, and with slight abuse of notation, it also holds that
	\begin{align*}
		& W_{U}^{\T} \cdot \left(n^{1/2}U_{i}\right)\\
		&\hspace{1em}= \left(R_{\wh{U}} R_{\wh{U}}^{\T} \right) \cdot \wh{U}^{\T} U \cdot \left(n^{1/2} U_{i}\right) + \ohp\left\{(n\rho)^{-1/2}\right\}\\
		&\hspace{1em}= R_{\wh{U}} \left(n^{1/2} R_{\wh{U}}^{\T} \wh{U}^{\T} \right) \cdot \left(U  W_{\Upsilon} W_{\Upsilon}^{\T}\right) \cdot U_{i} + \ohp\left\{(n\rho)^{-1/2}\right\}\\
		&\hspace{1em}= R_{\wh{U}} \left(n^{1/2} R_{\wh{U}}^{\T} \wh{U}^{\T} \right) \cdot \Xi \left(\Xi^{\T}\Xi\right)^{-1} \cdot \Xi_{i} + \ohp\left\{(n\rho)^{-1/2}\right\}\\
		&\hspace{1em}= R_{\wh{U}} \left(n^{1/2} R_{\wh{U}}^{\T} \wh{U}^{\T} \right) \cdot Z T_{\Xi}^{\T} \left(T_{\Xi} Z^{\T} Z T_{\Xi}^{\T}\right)^{-1} \cdot T_{\Xi} Z_{i} + \ohp\left\{(n\rho)^{-1/2}\right\}\\
		&\hspace{1em}= R_{\wh{U}} \left(n^{1/2} R_{\wh{U}}^{\T} \wh{U}^{\T} \right) \cdot Z \left(Z^{\T} Z\right)^{-1} \cdot Z_{i} + \ohp\left\{(n\rho)^{-1/2}\right\}\\
		&\hspace{1em}= R_{\wh{U}} \left(n^{-1} Z^{\T} \wh{Z}\right)^{\T} \cdot \left(n^{-1} Z^{\T} Z\right)^{-1} \cdot Z_{i} + \ohp\left\{(n\rho)^{-1/2}\right\} \\
		&\hspace{1em}= R_{\wh{U}} \Pi_{Z}^{\T} \cdot \breve{Z}_{i} + \ohp(1),
	\end{align*}
	where $\Pi_{Z} \cdot \left(n^{-1} Z^{\T} \wh{Z}\right)^{\T} \cdot Z_{i}$ converges to $\diag(\pi)^{1/2}Z_{i}$ given $Z_{i}=e_{\ell}$ by an analogous argument as in the preceding proof, leveraging the asymptotic normality arising from \cref{eq:prf_biom_vector_extension_directed}. Consequently, there exists a sequence of signed permutation matrices $(\Pi_{Z}) \subset \signorthGroup(k)$ such that for each fixed index $i$, the sequence of random vectors
	\begin{equation*}
		(n\rho)^{1/2}\left\{\Pi_{Z} \cdot \wh{Z}_{i} - \breve{Z}_{i} \mid Z_{i} = e_{\ell} \right\}
	\end{equation*}
	converges in distribution to a multivariate Gaussian random vector with mean zero and covariance matrix given by
	\begin{equation*}
		\left\{\diag(\pi_{Z})^{-1/2} \cdot T_{\Xi}^{-1} \cdot \Delta_{\Upsilon}^{-1}\right\} \cdot \operatorname{E}_{\upsilon}\left\{ g_{\rho_{\infty}}(\upsilon, T_{\Xi} e_{\ell}) \right\} \cdot \left\{\Delta_{\Upsilon}^{-1} \cdot (T_{\Xi}^{-1})^{\T} \cdot \diag(\pi_{Z})^{-1/2}\right\}.
	\end{equation*}
	By transposing the proof arguments to instead analyze $\wh{Y}_{j}$, there exists a sequence of signed permutation matrices $(\Pi_{Y}) \subset \signorthGroup(k)$ such that for each fixed index $j$, the sequence of random vectors of the form
	\begin{equation*}
		(n\rho)^{1/2}\left\{\Pi_{Y} \cdot \wh{Y}_{j} - \breve{Y}_{j} \mid Y_{j} = e_{\ell} \right\}
	\end{equation*}
	converges in distribution to a multivariate Gaussian random vector with mean zero and covariance matrix given by
	\begin{equation}
		\label{eq:prf_covariance_DirSBM_Yhat}
		\left\{\diag(\pi_{Y})^{-1/2} \cdot T_{\Upsilon}^{-1} \cdot \Delta_{\Xi}^{-1}\right\} \cdot \Ex_{\xi}\left\{ g_{\rho_{\infty}}(\xi, T_{\Upsilon} e_{\ell}) \right\} \cdot \left\{\Delta_{\Xi}^{-1} \cdot (T_{\Upsilon}^{-1})^{\T} \cdot \diag(\pi_{Y})^{-1/2}\right\}.
	\end{equation}
	To reiterate, here $\pi_{Z}$ and $\pi_{Y}$ denote the block membership probability vectors for the left and right latent factor vectors $Z_{i}$ and $Y_{j}$. The matrices of coefficients $T_{\Xi}$ and $T_{\Upsilon}$ defined in terms of $B$ satisfy $B = T_{\Xi}^{\T}T_{\Upsilon}$. Here, $\Delta_{\Xi} = \Ex(\xi\xi^{\T})$ where $\xi \sim T_{\Xi}\mathcal{G}$ and $\Delta_{\Upsilon} = \Ex(\upsilon\upsilon^{\T})$ with $\upsilon \sim T_{\Upsilon}\mathcal{H}$, while $g_{\rho_{\infty}}(\alpha, \beta) = (\beta^{\T}\alpha) \cdot (1-\rho_{\infty}\beta^{\T}\alpha) \cdot \alpha\alpha^{\T}$ depends on $\rho_{\infty} \in \{0,1\}$ and takes $k$-dimensional input vectors $\alpha, \beta$. This concludes the proof.
\end{proof}

\subsection{Asymptotic normality in varimax-transformed degree-corrected graph embeddings}
\label{sec:proof_vmx_dcsbm}
\begin{proof}[outline of Theorem~3]
	The proof proceeds in the same manner as in \cref{sec:proof_vmx_undirected_sbm} but instead by replacing $Z_{i}$ with $\theta_{i} \cdot Z_{i}$, properly scaling, and then finally integrating over all values $b \in \operatorname{supp}\mathcal{F}$. The compactness assumption precludes potential degenerate behavior near the boundary. Here, the integration step arises due to the general inability to condition on events of the form $\{\theta_{i} \cdot Z_{i}=x\}$ for fixed vectors $x$, since the distribution $\mathcal{F}$ need not be discrete.
\end{proof}


\end{document}